\documentclass[12pt]{l4dc2021} 
\usepackage{mathrsfs}
\usepackage{selectp}

\newcommand{\e}{\varepsilon}
\newcommand{\diam}{\textnormal{diam }}
\newcommand{\diag}{\textnormal{diag}}

\newcommand{\E}{\mathbf{E}}

\DeclareMathOperator{\spn}{span}

\DeclareMathOperator{\im}{im}

\DeclareMathOperator{\tr}{tr}
\DeclareMathOperator{\VEC}{\mathsf{vec}}
\DeclareMathOperator{\dop}{\mathsf{D}}
\DeclareMathOperator{\I}{\mathtt{I}}
\DeclareMathOperator{\J}{\mathtt{J}}
\DeclareMathOperator{\diff}{\mathtt{d}}

\title[Uninformative Optimal Policies in Adaptive LQR]{On Uninformative Optimal Policies in Adaptive LQR with Unknown B-Matrix}
\author{%
 \Name{Ingvar Ziemann} \Email{ziemann@kth.se}\\
 \Name{Henrik Sandberg} \Email{hsan@kth.se}\\
 \addr KTH Royal Institute of Technology%
}

\begin{document}
\maketitle

\begin{abstract}%   <- trailing '%' for backward compatibility of .sty file
This paper presents local  asymptotic minimax regret lower bounds for adaptive Linear Quadratic Regulators (LQR). We consider affinely parametrized $B$-matrices and known $A$-matrices and aim to understand when logarithmic regret is impossible even in the presence of structural side information. After defining the intrinsic notion of an \emph{uninformative optimal policy} in terms of a singularity condition for Fisher information we obtain local minimax regret lower bounds for such uninformative instances of LQR by appealing to van Trees' inequality (Bayesian Cram{\'e}r-Rao) and a representation of regret in terms of a quadratic form (Bellman error). It is shown that if the parametrization induces an uninformative optimal policy, logarithmic regret is impossible and the rate is at least order square root in the time horizon. We explicitly characterize the notion of an \emph{uninformative optimal policy} in terms of the nullspaces of system-theoretic quantities and the particular instance parametrization.
%In particular, if the entire matrix $B$ is unknown, the optimal policy can be uninformative if and only if the dimension of the input space exceeds the dimension of the state space, suggesting that regret of order $\sqrt{T}$ is essentially a high-dimensional phenomenon.
\end{abstract}

\section{Introduction}
Possibly given some structural side information, what is the asymptotic order of magnitude of regret for a fixed unknown instance of the linear quadratic regular (LQR)? We introduce a framework for thinking about this question in terms of the Fisher information matrix about an underlying parameter, $\theta$, generated by the optimal policy. This information quantity captures how much we learn about the underlying parameter by playing optimally in each round.   It is shown that when this information matrix, depending only on the optimal policy, satisfies a certain degeneracy property, logarithmic regret is impossible and the order of magnitude must be square-root in the observation horizon. Our argument relies on information comparison; we observe that a low regret policy must yield information about $\theta$ comparable to the optimal policy. That is, we use Fisher information to capture the \emph{exploration-exploitation} trade-off in adaptive LQR in terms of regret lower bounds.

Moreover, this reliance on Fisher information allows us to give conditions for when structural side-information is insufficient for logarithmic regret. If accurate structural  models are available to describe certain systems, one asks what the impact of such structure -- side information -- may be on the efficiency of learning algorithms. For instance, if we are controlling a networked system, we may want to impose graph structure. One can also imagine physical constraints imposing some symmetry or relation between certain coefficients.

%More broadly, the practical relevance of the above question arises in designing model-based reinforcement learning algorithms.  For us, the linear-quadratic setting serves as a theoretically tractable example for investigating this and studying reinforcement learning in continuous state and action spaces  \citep{recht2019tour, matni2019self}. 

Indeed, depending on the problem structure, it has been observed that there is a sharp regret phase-transition in learning linear quadratic regulators in its asymptotic scaling with the time horizon, $T$. There are several upper bound results in the $\sqrt{T}$ regime of regret when the entire system is unknown \citep{abbasi2011regret, mania2019certainty, cohen2019learning}, with \cite{simchowitz2020naive} establishing also an up-to-logs matching lower bound. However, in the presence of certain side-information or structure, logarithmic rates are attainable \citep{faradonbeh2020adaptive, cassel2020logarithmic}. The question of precisely specifying what structure makes logarithmic regret possible however remains open. We attempt to adress this by showing when it is not.

\paragraph{Contribution.}

We provide a natural framework for capturing the phase transition in regret depending on the structure of a nominal instance of LQR and possible side information encoded by an affine map.  When the Fisher information of the trajectory corresponding to optimal (stationary) policy for this instance satisfies a certain singularity condition, we say that the instance is \emph{uninformative} (see Section~\ref{uninfsec}). Assuming that the $A$-matrix is known and that $\VEC B = L\theta+B_0$, we establish that uninformativeness is a sufficient condition for logarithmic regret to be impossible. Indeed, our Theorem~\ref{thetheorem} shows that for all uninformative instances, regret $R_T$ (see (\ref{regretdef})), satisfies for some constant $C(\theta,\e)>0$
\begin{align}
\label{theprizehen}
\limsup_{T\to \infty}  \sup_{\theta' \in B(\theta,\e T^{-1/4})} \frac{R_T^\pi(\theta')}{\sqrt{T}} \geq C(\theta,\e).
\end{align}

To establish this result, a number of intermediate observations are made. First, we extend the exact (modulo terms $O(1)$) representation of regret as the cumulative Bellman error to the LQR setting in Lemma~\ref{testqlemma}. We then interpret this Bellman error, which is a sum of quadratic forms, as a sequence of estimation variances. The minimax regret, the left hand side of (\ref{theprizehen}), is then lower bounded by placing a suitable family of priors over shrinking subsets of $B(\theta,\e)$, to which we then apply Van Trees' inequality (Bayesian Cram{\'e}r-Rao). This interpretation is formalized in Lemma~\ref{decouplemma}, which lower bounds (\ref{theprizehen}) by inverse Fisher information. The final observation is that Fisher information itself has a relationship to regret in LQR. Specifically, Lemma~\ref{informationcomparisonlemma} establishes a principle of information comparison by spectral perturbation \citep{davis1970rotation, wedin1972perturbation, cai2018rate}. We show that any policy which has low regret automatically has Fisher information comparable to the Fisher information corresponding to use of the optimal policy (which is singular if the instance is uninformative). Roughly speaking, for any policy $\pi$, this idea states that
\begin{align*}
\textnormal{Fisher Information}(\pi) = \textnormal{Fisher Information}(\pi^*) + O(\textnormal{Regret})
\end{align*}
where $\pi^*$ is the optimal policy. Fisher information thus allows us to make precise the exploration-exploitation trade-off in adaptive LQR. In other words, the learner has to combat the singularity in Fisher information by adding excitation, however, since regret bounds information, this extra excitation necessarily has non-trivial contribution to regret. We believe this gives an appealing framework for understanding regret in online control.

\subsection{Related Work}
The problem of adaptively controlling an unknown instance of a linear quadratic system has a rather long history and dates back to at least \cite{aastrom1973self}. Early works \citep{goodwin1981discrete, lai1982least, campi1998adaptive} only asked that the adaptive algorithm be asymptotically optimal on average, which in our modern language can be rephrased as having sublinear regret. Historically, in this setting the emphasis on regret minimization as a means to study the performance of an adaptive algorithm appears first in the works by \cite{lai1986asymptotically} and \cite{lai1986extended}. See also \cite{guo1995convergence}. Noteably, in these last three works regret is shown to scale logarithmically with time, albeit subject to rather strong structural conditions. 

These logarithmic rates are of course in stark contrast with the current trend, where the emphasis has been on establishing bounds in the regime $\sqrt{T}$ (and with high probability). The present incarnation of this problem, the adaptive LQR model, was essentially popularized by \cite{abbasi2011regret} in which the authors produced an algorithm with $\tilde O(\sqrt{T})$ regret. A line of work following that publication focuses on improving and providing more computationally tractable algorithms in this setting \citep{ouyang2017control, dean2018regret, abeille2018improved, abbasi2019model, mania2019certainty, cohen2019learning, faradonbeh2020input,abeille2020efficient}. The emphasis of these works is entirely on providing upper bounds. Recently, some effort has been made to understand the complexity of the problem in terms of lower bounds. In particular,  \cite{simchowitz2020naive} provides matching (modulo constants and log factors) upper and lower bounds scaling correctly with the dimensional dependence given that the entire set of parameters $(A,B)$ are unknown. The works of \cite{cassel2020logarithmic} and \cite{ziemann2020phase} are also interesting in that they, for certain specific cases, provide $\sqrt{T}$ lower bounds that take the structure of the problem into account. Further, \cite{cassel2020logarithmic} shows that when the $A$-matrix is known and the optimal policy satisfies a certain non-degeneracy condition logarithmic rates are in fact achievable and a similar observation is made by \cite{ziemann2020phase} in the context of memoryless systems.

Our lower bound shares some features with that in \cite{simchowitz2020naive}, such as dimensional dependence, which also considers a local minimax version of regret. However, their bounds do not apply to the situation in which the learner is presented with structural side-information. Our proof approach is different and relies on Van Trees' inequality \citep{van2004detection,bobrovsky1987some}. This necessarily involves the Fisher information, which, quite naturally, allows for taking problem structure into account by considering different parametrizations of the problem dynamics. Cram{\'e}r-Rao type bounds have previously been used to establish regret lower bounds in adaptive LQR \citep{ziemann2020phase, ziemann2020regret}. However, these papers consider either restricted (memoryless) structure or do not take into account the possible singularity of Fisher information, which very much drives our result. We also note that the idea to bound a minimax complexity by a suitable family of Bayesian problems is well-known in the statistics literature \citep{gill1995applications}. See also \cite{van2000asymptotic,tsybakov2008introduction, ibragimov2013statistical} and the references therein.

Indeed, the adaptive control problem is intimitely connected to parameter estimation \citep{polderman1986necessity}. From the outset algorithm design has to a large extent been based on certainty equivalence; that is, estimating the parameters and plugging these estimates into an optimality equation, as if they were the ground truth \citep{aastrom1973self, mania2019certainty}. Our lower bound condition, \emph{uninformativeness}, is related to identifiability and inspired by a similar phenomenon in point estimation, which may become arbitrarily hard when Fisher information is singular \citep{rothenberg1971identification, goodrich1979necessary, stoica1982non, stoica2001parameter}. We also note that non-singularity of Fisher information is strongly related to the size of the smallest singular value of the covariates matrix in linear system identification \citep{faradonbeh2018finite, pmlr-v75-simchowitz18a,sarkar2019near,jedra2020finite}, which actually quantifies the corresponding rate of convergence \citep{jedra2019sample}.

Low regret linear quadratic control also fits into a wider context of online decision-making. \cite{lai1985asymptotically, burnetas1996optimal} solve asymptotically a set of problems known as bandits. Similar to LQR, depending on the problem structure and regret definition, the complexity of these problems also switches between $\log T$ and $\sqrt{T}$ \citep{shamir2013complexity}. Bandits are in some sense memoryless Markov Decision Processes (MDP). For finite state and action spaces, \cite{burnetas1997optimal}, characterizes instance specific regret for MDPs. See also \cite{ok2018exploration} and the referencs therein. Returning to the control context, there have also been recent advances in the partially observed adaptive Linear Quadratic Gaussian (LQG) setting, \citep{lale2020logarithmic} and unknown cost LQR with adversarial disturbances \citep{hazan2020nonstochastic}. We note that logarithmic rates reappear in the partially observed case due to certain excitation properties of the output sequence. Moreover, there is an interesting parallel line of work which emphasizes robustness in adaptive LQR \citep{dean2019sample, umenberger2019robust}. See also \cite{recht2019tour, matni2019self} and the references therein for an overview of the relationship between control and learning.

\subsection{Outline}
We first define the problem in Section~\ref{probfor}. Section~\ref{regsec} relates regret to sub-optimal solutions of the Bellman equation. Section~\ref{infosec} reviews Fisher information and studies \emph{uninformativeness}, which is the key notion used in Section~\ref{mainsec} to derive our regret lower bound. Finally Section~\ref{consec} concludes. All proofs can be found in the arXiv version of this paper \citep{ziemann2020uninformative}, which includes our appendices.

\paragraph{Notation.} 
We use $\succeq$ (and $\succ$) for (strict) inequality in the matrix positive definite partial order. By $\|\cdot\|$ we denote the standard $2$-norm by $\|\cdot\|_\infty$ the matrix operator norm (induced $l^2 \to l^2$) and $\rho(\cdot)$ denotes the spectral radius. Moreover, $\otimes$, $\VEC$ and $\dagger$ are used to denote the Kronecker product, vectorization (mapping a matrix into a column vector by stacking its columns), and the Moore-Penrose pseudoinverse, respectively. For a sequence of vectors $\{v_t\}_{t=1}^n, v_t \in \mathbb{R}^d$ we use $v^n=(v_1,\dots,v_n)$ defined on the $n$-fold product $\mathbb{R}^{d\times n}$. The set of $k$-times continuously differentiable functions on $\mathbb{R}^{d}$ is denoted $C^k(\mathbb{R}^d)$. We use $\dop$ for Jacobian, $\diff$ for differential and $\nabla$ for the gradient. That is, for a scalar function, $\nabla f$ denotes a column vector of its first derivatives. We write $\E$ for the expectation operator, with superscripts indicating policy, and subscripts indicating parameters.

\section{Problem Formulation}
\label{probfor}

Fix an unknown parameter $\theta \in \Theta$ where $\Theta$ is an open subset of $\mathbb{R}^{d_\theta}$. Let $\{x_t\}$ be a controlled process on $ \mathbb{R}^{d_x}$ with dynamics
\begin{align}
\label{themodel}
x_{t+1} &= Ax_t+B(\theta)u_t+w_t, & x_0&=0, & t&=0,1,2,\dots
\end{align}
with control process $\{u_t\}$ on $ \mathbb{R}^{d_u}$, additive noise process $\{w_t\}$ on  $\mathbb{R}^{d_x}$ so that $A \in \mathbb{R}^{d_x\times d_x}$ and $B=B(\theta)\in \mathbb{R}^{d_x\times d_u}$. It is assumed that $B$ depends affinely on $\theta$, to be made precise momentarily. We consider the stage cost function $c(x,u) = x^\top Q x+u^\top R u$, with $Q \in \mathbb{R}^{d_x\times d_x}, R \in \mathbb{R}^{d_u\times d_u}$. We further denote the $\sigma$-field generated by $x_1,\dots,x_t$ and possible auxilliary randomization by $\mathcal{F}_t$.  The adaptive control objective is to design a policy $\pi$ that minimizes the cumulative cost
\begin{align}
\label{costv}
V^\pi_T(\theta) = \sum_{t=0}^{T-1} \E^\pi_{\theta} c(x_t,u_t) =  \sum_{t=0}^{T-1} \E^\pi_{\theta} \left[  x_t^\top Q x_t+u_t^\top R u_t \right]
\end{align}
without a priori knowledge of the parameter $\theta$. This paper studies the fundamental limitations to this problem. We make the following standing assumptions about (\ref{themodel})-(\ref{costv}).

\begin{enumerate}
\item[A1.] $\{w_t\}$ is iid $p(\cdot)$ with $\E w_t w_t^\top =\Sigma_w$ and $p(\cdot) \in C^1(\mathbb{R}^{d_x})$ with finite Fisher information.
\item[A2.]  The cost function is strongly convex. More precisely, $Q \succ 0$ and $R\succ 0$.
\item[A3.] The dynamics (\ref{themodel}) are stabilizable at $\theta$; there exists $K(\theta)$ with $\rho(A-B(\theta)K(\theta))< 1$.
\item[A4.] $B(\theta)$ is affine in $\theta$; $\VEC B(\theta) = L\theta +\VEC B_0$ for  matrices $L\in\mathbb{R}^{d_xd_u \times  d_\theta}, B_0 \in \mathbb{R}^{d_x\times d_u}$.
\end{enumerate}

We say that a tuple $(\theta, A, B(\cdot), Q,R,p)$ satisfying the above assumptions is a parametrized instance of LQR. Our goal will be to devise regret lower bounds in terms of the nominal instance, $\theta$, and the structure of the parametrization, i.e, in terms of the matrix $L$. To motivate why this flexibility in $L$ is beneficial, consider the following example.

\begin{example}
\label{ch6:graphex}
Consider the system (\ref{themodel}) and suppose that it is known that $B$ has network strucutre, say, it is the Laplacian of some known graph $\mathtt{G}$ with unknown weights. For simpicity let us assume that $d_x = d_u = 3$ and that the graph adjacency matrix is given by
\begin{align*}
\mathtt{Adj} = \begin{bmatrix}
0 & 1 & 1\\
1 & 0 & 0\\
1 & 0 &0
\end{bmatrix}
\end{align*}
so that nodes 1 and 2 and 1 and 3 are connected by an edge. We assume that this topological information is available to the user. The Laplacian is thus of the form
\begin{align*}
B = \begin{bmatrix}
\theta_1 & \theta_2  & \theta_3 \\
\theta_2 & \theta_4 & 0\\
\theta_3 & 0 & \theta_5
\end{bmatrix} \textnormal{ with }
\VEC B = \begin{bmatrix}
\theta_1 & \theta_2 & \theta_3 & \theta_2 &\theta_4 & 0 &\theta_3 & 0 &\theta_5
\end{bmatrix}^\top.
\end{align*}
Clearly, $\VEC B = L \theta$, where $\theta \in \mathbb{R}^5$ for some matrix $L\in \mathbb{R}^{9\times 5}$. From an identification perspective, one thus suspects that the side-information of knowing the network topology reduces the hardness of the adaptive control task. Namely, the number of parameters to be identified is only $d_\theta = 5$ instead of $9 = \textnormal{number of entries}(B)$.
\end{example}

\paragraph{Optimality and Regret.}
The objective (\ref{costv}) is equivalent to minimizing the regret
\begin{align}
\label{regretdef}
R_T^\pi(\theta) = V_T^\pi(\theta)-V_T^*(\theta)
\end{align}
where $V_T^*(\theta) = V_T^{\pi^*}(\theta)$ is the cumulative cost corresponding to the optimal policy $\pi^*(\theta)$, defined implicitly by the Riccati recursion \citep{bertsekas1995dynamic}.

%\begin{remark}
%Certain authors refer to the quantity (\ref{regretdef}) as the expected regret. However, as all our expressions are under the expectation, we shall for brevity refer to (\ref{regretdef}) as simply the regret.
%\end{remark}

\section{Regret and Cumulative Bellman Error}
\label{regsec}
One expects that $V_n^*(\theta) \approx  V_n^{\pi^\infty}(\theta)$ where $\pi^\infty$ instead of relying on the Riccati recursion relies on its stationary limit. The stationary policy $\pi^\infty$ is then characterized by minimization of the quadratic form
\begin{align*}
\phi(x,u, \theta) = x^\top Q x + u^\top R u +  [A(\theta)x+B(\theta)u]P(\theta)[A(\theta)x+B(\theta)u)]-x^\top P(\theta) x.
\end{align*}
If $P(\theta)$ solves the discrete algebraic Riccati equation $\min_{u\in \mathbb{R}^{d_u}} \phi(x,u,\theta) =0$ is equivalent to the average cost Bellman equation for this problem \citep{bertsekas1995dynamic}. Our first theorem exploits the quadratic nature of $\phi$ to transform the problem essentially into one of sequential estimation. 
\begin{theorem}
\label{regexpthm}
Under assumptions A1-A4 it holds that
\begin{multline}
\label{ziemannregexp}
R_T^\pi(\theta) =\sum_{t=0}^{T-1} \E_{\theta}^\pi\left(\phi(x_t,u_t, \theta)\right)+O(1)\\
=\sum_{t=0}^{T-1} \E_{\theta}^\pi\left( (u_t-K(\theta)x_t)^\top\big[R +B^\top(\theta) P(\theta) B(\theta)\big](u_t-K(\theta)x_t)\right)+O(1)
\end{multline}
where the term $O(1)$ is uniformly bounded in the stability region of the optimal stationary policy $K(\theta)$, i.e. on each nonempty set $\{\theta' \in \mathbb{R}^{d_\theta} | \rho(A-B(\theta')K(\theta)) \leq 1 -\zeta\}$, $\zeta \in (0,1)$.
\end{theorem}
See Appendix~\ref{bkappendix} for the proof. Note that such a neighborhood always exists for some $\zeta$ by upper semi-continuity of the spectral radius and stability of $A-B(\theta)K(\theta)$. The proof idea behind Theorem~\ref{regexpthm} is due to  \cite{burnetas1997optimal}.  Similar expressions are also obtained in e.g. \cite{faradonbeh2020adaptive, simchowitz2020naive}. 

 %The essential ingredient is the exponential convergence of the Riccati recursion to its stationary limit, see e.g. \cite{anderson2012optimal} which relates $R_T^\pi$ to a sum over $\phi$. 

Notice that when the choice of $u_t$ is made $x_t$ is known. Hence (\ref{ziemannregexp}) informs us that regret minimization is essentially equivalent to minimizing a cumulative weighted estimation error for the sequence of estimands $K(\theta)x_t$. To be clear, our perspective is that we wish to estimate the function value of $\theta \mapsto K(\theta)x_t$ where the function $K(\theta)x_t$ is revealed at time $t$ by virtue of observations of the $x_t$. Moving to a Bayesian setting, the entire trajectory $(x^{T+1},u^T)$ is then interepreted as a noisy observation of the underlying parameter $\theta$. A natural approach for variance lower bounds is to rely on Fisher information and use Cram{\'e}r-Rao type bounds.

\section{Fisher Information Theory}
\label{infosec}

Let us recall the definition of Fisher information. For a parametrized family of probability densities $\{q_\theta, \theta \in \Theta\}$, $\Theta \subset \mathbb{R}^{d}$, Fisher information $\I_p(\theta)\in \mathbb{R}^{d\times d}$ is 
\begin{align}
\label{fisherdef}
\I_q(\theta)= \int \nabla_\theta \log q_\theta(x)\left[\nabla_\theta \log q_\theta(x)\right]^\top q_\theta(x) dx
\end{align}
whenever the integral exists. For a density $\lambda$, we also define the location integral
\begin{align}
\label{fisherlocdef}
\J(\lambda)= \int \nabla_\theta \log \lambda(\theta) \left[\nabla_\theta \log \lambda(\theta) \right]^\top \lambda(\theta) d\theta.
\end{align}
Again, provided of course that the integral exists. See \cite{ibragimov2013statistical} for details about these integrals and their existence.

We now study Fisher information where $q$ in (\ref{fisherdef}) is the joint density of $x^{T+1}$.
\begin{lemma}
Under Assumption 1, Fisher information about $\theta$ given observations of $(x^{T+1},u^T)$ in the model (\ref{themodel}) is given by
\begin{align}
\label{problemspecificinformation}
 \I^T(\theta; u^T) =\E \sum_{t=0}^T [\dop_\theta[B(\theta)u_t]] ^\top \J(p) \dop_\theta[B(\theta)u_t].
 \end{align}
\end{lemma}

See appendix \ref{chainfisherapp} for the proof. We now turn to investigating Fisher information corresponding to observations generated by the optimal policy $u_t = K(\theta)x_t$. It will be especially interesting to study when (\ref{problemspecificinformation}) becomes singular.

\subsection{Uninformative Optimal Policies}
\label{uninfsec}
Naively, the perspective discussed in Section~\ref{regsec} viewing (\ref{ziemannregexp}) as a cumulative estimation error suggests a lower bound on the scale $\log T$, since one might think that the errors variances should decay as $1/t$. However, when Fisher information is singular,  this reasoning might fail.

We say that an instance $(\theta, A, B(\cdot), Q,R,p)$ is \emph{uninformative} if 
\begin{enumerate}
\item $  \I_*^T(\theta) =\E \sum_{t=0}^T [\dop_\theta[B(\theta)K(\theta)x_t]] ^\top \J(p) \dop_\theta[B(\theta)K(\theta)x_t]$ is singular for all $T$ under closed loop dynamics
\begin{align*}
x_{t+1} = (A-B(\theta)K(\theta))x_t+w_t; \textnormal{ and}
\end{align*}
\item There exists a vector $\tilde v$ in the nullspace of $\I^T(\theta)$ such that $\dop_\theta \VEC K(\theta) \tilde v \neq 0$.
\end{enumerate}
In other words, uninformativeness stipulates that observation of an optimally regulated example is not sufficient to (locally) identify the optimal policy.

\paragraph{Algebraic Charaterization of Uninformativeness.}
The above description of what constitutes an uninformative optimal policy is somewhat indirect. We can characterize uninformativeness directly in terms of the instance parameters. 

\begin{proposition}
\label{charprop}
The instance  $(\theta, A, B(\cdot), Q,R,p)$ is uninformative if and only if there exists a vector $\tilde v$ such that
\begin{align}
\label{uninformativechar}
\begin{cases}
\tilde v &\in \ker L^\top [K K^\top  \otimes  \J(p) ]L\\
\tilde v&\notin \ker \dop_\theta \VEC K(\theta)
\end{cases}
\end{align}
where $L$ is such that $\VEC B = L\theta+B_0$.
\end{proposition}
The proof is given in Appendix~\ref{charpropproof}. Any subspace $\mathtt{U}$ of maximal dimension for which all nonzero $\tilde v \in \mathtt U$ satisfy (\ref{uninformativechar}) is called an \emph{information singular subspace}. We will later see that the dimension of any such subspace gives the dimensional dependence in our lower bound.

As \cite{cassel2020logarithmic} show that logarithmic regret rates for known $A$-matrix and unknown $B$-matrix are attainable if $KK^\top$ is nonsingular, it is interesting to note that our notion of uninformativeness immediately rules out the possibility of nonsingular $KK^\top$.

\begin{corollary}
\label{casselcorr}
If the instance $(\theta, A, B(\cdot), Q,R,p)$ is uninformative with $\VEC B = \theta$ so that $L=I_{d_\theta}$, then $KK^\top$ is singular.
\end{corollary}

The following examples illustrate the concept in terms of certain simple model structures.

\begin{example}
\label{scalarexample}
Consider a ``scalar'' LQR, with nonzero $A=a \in \mathbb{R}$ known, and $B=\theta \in \mathbb{R}$ unknown. Since the optimal linear feedback law is $0$ if and only if $\theta =0$, it follows that scalar LQR is uninformative if and only if the input matrix is $B=\theta=0$. Notice that scalar $B \approx 0$ is precisely the construction used in the lower-bound proof of \cite{cassel2020logarithmic}.
\end{example}

In the next example, we illustrate that one can explicitly compute the dimension of $\mathtt{U}$, the number of parameters not excited by the optimal policy.

\begin{example}
\label{memorylessexample}
Consider a ``memoryless'' linear quadratic regulator (c.f. \cite{ziemann2020phase})
\begin{align}
\label{memorylessdynamics}
\begin{bmatrix}
r_{t+1}\\
y_{t+1}
\end{bmatrix}
=
\begin{bmatrix}
G & 0 \\
I & 0
\end{bmatrix}
\begin{bmatrix}
r_{t}\\
y_{t}
\end{bmatrix}
+
\begin{bmatrix}
0 \\
-F
\end{bmatrix}
u_t
+\begin{bmatrix}
n_{t}\\
v_{t}
\end{bmatrix},
\end{align}
with $R=I_{d_u}, Q=I_{d_x}$. Due to the memoryless property in the second $x$-coordinate $y_t$, the optimal regulation of this instance can be reduced to
\begin{align}
\label{memorylessexamplemodel}
\begin{cases}
y_{t+1} =r_t- Fu_t +v_t\\
\min_{\{u_t\}} \E \sum_{t=1}^T  \|r_t -Fu_t\|^2 +  \|u_t\|^2
\end{cases}
\end{align}
which has optimal policy in feedforward form $u_t = (F^\top  F+I_{d_u})^{-1} F^\top  r_t$. Suppose that $\VEC F = \theta$, so that the structure (\ref{memorylessdynamics}) is known. It can be shown this instance is uninformative if and only if $KK^\top$ is singular. Moreover, $\dim \mathtt{U} =d_y \times \dim\ker KK^\top $. See the appendix for the proof.
\end{example}

Note that in this case, the condition agrees exactly with that in \cite{cassel2020logarithmic}. 
\section{Main Result}
\label{mainsec}
Our main result is a \emph{local} asymptotic minimax regret lower bound. It considers a shrinking neighborhood in parameter space around $\theta$ and states that if the nominal instance $\theta$ is \emph{uninformative} any policy $\pi$ suffers regret on the order of magnitude $\sqrt{T}$ on at least one instance of this neighborhood.

\begin{theorem}
\label{thetheorem}
Assume that A1-A4 hold and suppose that the nominal instance \\ $(\theta, A, B(\cdot), Q,R,p)$ is uninformative. Then for any $\Gamma\succ 0$ and any $\e>0$ such that $B(\theta,\e) \subset \{\theta' \in \mathbb{R}^{d_\theta} |\rho(A-B(\theta)K(\theta')) <1 \}$, any policy $\pi$ satisfies
\begin{align}
\label{regretlowerblin}
\limsup_{T\to \infty}  \sup_{\theta' \in B(\theta,\e T^{-1/4})} \frac{R_T^\pi(\theta')}{\sqrt{T}} \geq C(\theta,\e,\Gamma)
\end{align}
where
\begin{equation}
\begin{aligned}
\label{optconstantC}
C(\theta,\e,\Gamma)=\inf_{C' > 0} \max\Bigg\{p^2 \tr \Bigg(\frac{ \Gamma}{F(\theta,\e,C')} \otimes\big[R +B^\top( \theta) P( \theta) B(\theta)\big]\\
\times {\dop_\theta\VEC K(\theta) \tilde W_0 \tilde W_0^\top[\dop_\theta\VEC K(\theta)]^\top}
\Bigg), C'\Bigg\}
\end{aligned}
\end{equation}
where $p = \liminf \mathbf{P}^\pi_{\theta'\sim \lambda}\left(\left\{\sum_{t=k\lceil\sqrt{T}\rceil}^{(k+1)\lfloor\sqrt{T}\rfloor}x_t x_t^\top \succeq \sqrt{T} \Gamma\right \}\right) $, and where  $\tilde W_0$ is an orthonormal matrix with columns spanning an eigenspace of dimension $\lceil\dim \mathtt{U}/2\rceil$ of $\mathtt{U}$, and further
\begin{align*}
F(\theta,\e,C')&=\|L\|_\infty^2  \|Q^{-1}\|_\infty g^{\pi^*}(\theta)   \frac{ \tr \J(p)}{\lceil\dim \mathtt{U}/2\rceil}\|\dop^2_\theta \tr K(\theta)K^\top(\theta) \|_\infty\times \frac{\e^2}{2}\\
&+2 \| L\|^2_\infty \| [R+B^\top(\theta) P(\theta) B(\theta)]^{-1}\|_\infty \frac{ \tr \J(p)}{\lceil\dim \mathtt{U}/2\rceil} C'+\|\J(\lambda)\|_\infty
\end{align*}
where the infimum with respect to $\pi$ is taken subject to $\limsup_{T\to \infty}  \frac{R_T^\pi(\theta)}{\sqrt{T}} \leq C'$ and $F$ is defined for any $\lambda \in C_c^\infty[B(\theta,\e) ]$.
%where $\|\J(\lambda(1,1,2))\|_\infty$ is a fixed constant (see Appendix~\ref{molliapp}).
\end{theorem}

The bound can be optimized in terms of the state covariance proxy $\Gamma$. For instance, the choice $\Gamma=\Sigma_w$ results in $p=1$ in (\ref{optconstantC}). The quotient $\tr \Gamma \big/ F(\theta,\e,C_\pi)$ is approximately the trace of the state (observation) variance divided by an upper bound for inverse normalized Fisher information and as such can be interpreted as a signal-to-noise ratio (SNR). The optimization problem (\ref{optconstantC}) can thus be understood as to balance high SNR with good control performance. Indeed, this balancing constitutes one of the key observations leading to the proof of Theorem~\ref{thetheorem} (Appendix~\ref{mainthmproofappendix}); if regret is $O(C_\pi \sqrt{T})$, the dominant component of $F(\theta,\e,C_\pi)$ is $O(C_\pi)$. 

We also note in passing that together with the scalar case, Example~\ref{scalarexample},  we immediately see that the global minimax regret, also with unknown $A$-matrix, is at least of order $\sqrt{T}$.

\paragraph{Dimensional Depedence.}
Observe that the optimization problem (\ref{optconstantC}) has value asymptotically proportional to 
\begin{multline}
\sqrt{\tr \Bigg(\underbrace{\left[\Gamma \otimes \big[R +B^\top( \theta) P(\theta) B( \theta)\big] \right]}_{\succ 0 } {\dop_\theta\VEC K(\theta) \tilde W_0 \tilde W_0^\top[\dop_\theta\VEC K(\theta)]^\top}\Bigg)}\Bigg/\sqrt{\frac{\tr \J(p)}{\dim \mathtt{U}}} \\
\propto \sqrt{\tr \Bigg({\dop_\theta\VEC K(\theta) \tilde W_0 \tilde W_0^\top[\dop_\theta\VEC K(\theta)]^\top}\Bigg)}\Bigg/\sqrt{\frac{\tr \J(p)}{\dim \mathtt{U}}}\\
\propto \sqrt{\frac{(\dim \mathtt{U})^2}{d_x}}
\end{multline}
where $\mathtt{U}$ is the information singular subspace (\ref{uninformativechar}). That is, the dimensional depedence is proportional to the size of the parameter subspace not excited by the optimal policy relevant for identification of $K(\theta)$. Notice further that if the entire matrix $B$ is unknown, $\mathtt{U}$ is a subspace of $\mathbb{R}^{d_x d_u}$ which further relates to the dimensions of the problem.

\paragraph{Outline of the Proof of Theorem~\ref{thetheorem}}
We provide a brief outline of the proof, of which the details are given in Appendix~\ref{mainthmproofappendix}. The following steps are key:

\begin{itemize}
\item We relate regret to a sequential estimation problem via cumulative Bellman error. This is the content of Theorem~\ref{regexpthm}.
\item We consider the $\sup$ of regret over a shrinking neighborhood. This $\sup$ is relaxed to a Bayesian estimation problem for an increasingly concentrated sequence of priors to which we apply van Trees' inequality.
\item The role of these increasingly concentrated priors is to allow us to consider the local properties of Fisher information on the domain of the prior.
\item A truncation argument is used to decouple conditional Fisher information from the state trajectory.
\item We establish a principle of information comparsion; small regret implies that the Fisher information of that policy is near singular under the hypothesis that the optimal policy is uninformative. This implies large estimation error and thus large regret. The proof of this principle relies on the Davas-Kahan $\sin \theta$ Theorem \citep{davis1970rotation}.
\end{itemize}

\section{Discussion and Conclusion}
\label{consec}
We have established a framework which is able to both qualitatively and quantitatively explain why $\sqrt{T}$ regret occurs in certain instances of adaptive linear quadratic control. Essentially, if the optimal policy does not sufficiently excite all directions of the underlying parameter space necessary for identification of that same optimal policy, logarithmic regret is impossible. We make this explicit by a singularity condition on that optimal policy's Fisher information which takes problem structure into account. Indeed, singular information has deep connections to classical notions of identifiability  \citep{rothenberg1971identification, goodrich1979necessary, stoica1982non, stoica2001parameter}, and we believe that further investigation of this relation would be an interesting line of future work. An extension to unknown $A$-matrix is of course also highly relevant. 

Moreover, the tightness of our notion remains an open question. Some preliminary evidence is provided by the results of \cite{cassel2020logarithmic} which establish logarithmic regret under the nonsingularity condition $KK^\top \succ 0$ in the case of known $A$-matrix. Our  Corollary~\ref{casselcorr} is consistent with their observations. This opens an exciting problem: to prove or disprove whether logarithmic regret is attainable if and only if the instance  is not uninformative.

\acks{This work was supported in part by the Swedish Research Council (grant 2016-00861), and the Swedish Foundation for Strategic Research (Project CLAS). The authors express their gratitude to Yishao Zhou for valuable input.}

\newpage

\bibliography{lqrlocminimax}

\newpage

\appendix
The appendix provides the proofs of all results stated in the main text. We also review some of the auxilliary material key to our results.

\tableofcontents

\newpage
\section{Proof of Theorem~\ref{thetheorem}}
Here we give the proof of Theorem~\ref{thetheorem} and two auxilliary lemmas. We work under the convention that assumptions A1-A4 apply.

\paragraph{Proof of Theorem~\ref{thetheorem}}
\label{mainthmproofappendix}
The first part of the proof consists of combining lemmas~\ref{decouplemma} and \ref{informationcomparisonlemma} to get an information inequality for regret. This inequality is difficult to evaluate directly, so we split the time horizon into $T$ into $\sqrt{T}$-many equal summands, on which information is nearly constant. To proceed with this argument, we now fix the scale of the estimation problem (\ref{td:1}) by rescaling our as prior $\lambda_T(\cdot) = \lambda(T^{1/4} \cdot)$. Observe that the prior information then scales as $\J(\lambda_T)= T^{1/2} \J(\lambda)$. For convenience we also denote $\e_T = \e T^{-1/4}$ which is the radius of support of the rescaled prior. Notice further that we may work under the hypothesis, for some constant $C'>0$, that
\begin{align}
\label{uniformlysmallregret}
\limsup_{T\to \infty}  \frac{R_T^\pi(\theta)}{\sqrt{T}} \leq C'
\end{align}
since otherwise the claim is obviously true.

\paragraph{Restricting to the Uninformative Space}
By virtue of Lemma~\ref{decouplemma}, we can relate regret to a sum over Fisher information as
\begin{multline*}
  \sup_{\theta' \in B(\theta,\e T^{-1/4})} R_T^\pi(\theta')\\
 \geq   \sum_{k=1}^{\lfloor\sqrt{T-1}\rfloor} [\mathbf{P}^\pi_{\theta'\sim \lambda}(A_k^T)]^2\tr \Bigg(\left(\Gamma \otimes \big[R +B^\top(\theta) P(\theta) B( \theta)+O(\e T^{-1/4} )\big]  \right) \left( \dop_\theta \VEC K( \theta)+ O(\e T^{-1/4})\right)\\
  \times 
   \left( \sqrt{T}\left( \E^\pi_{\theta'\sim \lambda} \I^{(k+1)\lfloor\sqrt{T}\rfloor}(\theta) + \J(\lambda)+\Psi \right)^{-1}  \right)\left(\dop_\theta \VEC K( \theta+O(\e T^{-1/4}))\right)^\top\Bigg)+O(1).
\end{multline*}
Spectrally decompose now
\begin{align}
\label{mainthmeq2}
\E^\pi_{\theta\sim \lambda_T}\left[ \I^{(k+1)\lfloor\sqrt{T}\rfloor}(\theta)\right]+\J(\lambda) +\Psi=  V_1(T,k) \Lambda_1(T,k) V_1^\top(T,k) + V_0(T,k) \Lambda_0(T,k) V_0^\top(T,k)
\end{align}
so that
\begin{multline*}
\left(\E^\pi_{\theta\sim \lambda_T}\left[ \I^{(k+1)\lfloor\sqrt{T}\rfloor}(\theta)\right]+\J(\lambda) +\Psi\right)^{-1}\\=  V_1(T,k) \Lambda_1^{-1}(T,k) V_1^\top(T,k) + V_0(T,k) \Lambda_0^{-1}(T,k) V_0^\top(T,k)
\end{multline*}
with the variables $V_i$ defined as in Lemma~\ref{informationcomparisonlemma}. That is, the columns of $V_0$ span the smallest eigenspaces of the left-hand side of (\ref{mainthmeq2}), and the columns of $V_1$ the largest. It follows that 
\begin{multline}
\label{combinethis}
  \sup_{\theta' \in B(\theta,\e T^{-1/4})} R_T^\pi(\theta') \geq  \sum_{k=1}^{\lfloor\sqrt{T-1}\rfloor} [\mathbf{P}^\pi_{\theta'\sim \lambda}(A^T_k)]^2\tr \Bigg(\left(\Gamma \otimes \big[R +B^\top(\theta) P(\theta) B( \theta)+O(\e T^{-1/4} )\big]  \right)\\
  \times \left( \dop_\theta \VEC K( \theta)+ O(\e T^{-1/4})\right)\\
  \times 
   \left( \sqrt{T}V_0(T,k) \Lambda_0^{-1}(T,k) V_0^\top(T,k) \right)\left(\dop_\theta \VEC K( \theta+O(\e T^{-1/4}))\right)^\top\Bigg)+O(1)
\end{multline}
Take $\Psi = \psi\tilde V_1 \tilde V_1^\top $ with $\psi \in \mathbb{R}$ and let $\psi \to \infty$ (the definition of $\tilde V_1$ is given in (\ref{specdecs})). Then since $\Psi$ was arbitrary, by  Lemma~\ref{informationcomparisonlemma} we find
\begin{multline*}
  \sup_{\theta' \in B(\theta,\e T^{-1/4})} R_T^\pi(\theta') \geq  \sum_{k=1}^{\lfloor\sqrt{T-1}\rfloor} [\mathbf{P}^\pi_{\theta'\sim \lambda}(A^T_k)]^2\tr \Bigg(\left(\Gamma \otimes \big[R +B^\top(\theta) P(\theta) B( \theta)+O(\e T^{-1/4} )\big]  \right)\\
  \times \left( \dop_\theta \VEC K( \theta)+ O(\e T^{-1/4})\right)\\
  \times 
   \left( \sqrt{T}\tilde V_0(T,k) \Lambda_0^{-1}(T,k) \tilde V_0^\top(T,k) \right)\left(\dop_\theta \VEC K( \theta+O(\e T^{-1/4}))\right)^\top\Bigg)+O(1)
\end{multline*}

Now, by Lemma~\ref{informationcomparisonlemma} the quantity $\sigma_{\lceil\dim \mathtt{U}/2\rceil}(\Lambda^{-1}_0(T,k))$ has lower bound
\begin{multline}
\label{sigmamininthm}
\liminf_{T\to \infty} \sqrt{T} \times  \sigma_{\lceil\dim \mathtt{U}/2\rceil}(\Lambda^{-1}_0((k+1)\lfloor\sqrt{T}\rfloor))\geq \frac{1}{F(\theta,\e,C')} \\
=
\lceil\dim \mathtt{U}/2\rceil\Bigg/ \Bigg\{ \Bigg(\|L\|_\infty^2  \|Q^{-1}\|_\infty g^{\pi^*}(\theta)  \tr [J(p)] \|\dop^2_\theta \tr K(\theta')K^\top(\theta') \|_\infty\times \frac{\e^2}{2}\\
+2 \| L\|^2_\infty \| [R+B^\top(\theta) K(\theta) B(\theta)]^{-1}\|_\infty \tr \left( \J(p)  \right)  C'+\|\J(\lambda)\|_\infty\Bigg)\Bigg\}
\end{multline}
since $\e_T = \e T^{-1/4}$, $\J(\lambda_T) = T^{1/2} \J(\lambda)$. Combining (\ref{combinethis}) and (\ref{sigmamininthm}), the above yields that
\begin{multline}
\label{finitetimebound}
 \sup_{\theta' \in B(\theta,\e T^{-1/4})} R_T^\pi(\theta')\geq \sum_{k=1}^{\lfloor\sqrt{T-1}\rfloor}\Bigg\{[\mathbf{P}^\pi_{\theta'\sim \lambda}(A^T_k)]^2 \tr \Bigg(\Bigg[\frac{\Gamma}{F(\theta,\e,C')}\\
 \otimes \big[R +B^\top(\theta) K(\theta) B( \theta)\big]\Bigg] {\dop_\theta\VEC K(\theta') \tilde W_0 \tilde W_0^\top[\dop_\theta\VEC K(\theta')]^\top} \Bigg) +o(1)\Bigg\}.
\end{multline}
Thus balancing terms, we either have large regret -- (\ref{uniformlysmallregret}) does not hold -- or that the term $F(\theta,\e,C')$ is very near zero. Wherefore, by evaluating the limit superior of (\ref{finitetimebound}),
\begin{align*}
\limsup_{T\to \infty}  \sup_{\theta' \in B(\theta,\e)} \frac{R_T^\pi(\theta')}{\sqrt{T}} \geq C(\theta,\e,\Gamma)
\end{align*}
with
\begin{align*}
C(\theta,\e,\Gamma)=\inf_{C' > 0} \max\Bigg\{p^2 \tr \Bigg(\frac{ \Gamma}{F(\theta,\e,C')} \otimes\big[R +B^\top( \theta) P( \theta) B(\theta)\big]\\
\times {\dop_\theta\VEC K(\theta) \tilde W_0 \tilde W_0^\top[\dop_\theta\VEC K(\theta)]^\top}
\Bigg), C'\Bigg\}
\end{align*}
where $p = \liminf \mathbf{P}^\pi_{\theta'\sim \lambda}(A^T_k) $. $\hfill\blacksquare$

\begin{lemma}
\label{decouplemma}
Fix $\e>0$, $\Gamma\succ 0$, and let  $\lambda \in C_c^\infty[B(\theta,\e) ]$, then
\begin{multline}
\label{timedecoup}
 \sup_{\theta' \in B(\theta,\e)} R_T^\pi(\theta')\geq   \sum_{k=1}^{\lfloor\sqrt{T-1}\rfloor} [\mathbf{P}^\pi_{\theta'\sim \lambda}(A_k)]^2\tr \Bigg(\left(\Gamma \otimes \big[R +B^\top(\theta) P( \theta) B(\theta)+O(\e)\big]  \right) \\
\times\left( \dop_\theta \VEC K( \theta)+O(\e) \right)  \left( 
 \sqrt{T}\left( \E^\pi_{\theta'\sim \lambda} \I^{(k+1)\lfloor\sqrt{T}\rfloor}(\theta) + \J(\lambda)+\Psi \right)^{-1}  \right)\left(\dop_\theta \VEC K(\theta)+O(\e)\right)^\top\Bigg)+O(1)
\end{multline}
where $A^T_k = \{\sum_{t=k\lceil\sqrt{T}\rceil}^{(k+1)\lfloor\sqrt{T}\rfloor}x_t x_t^\top \succeq \sqrt{T} \Gamma \}$.
\end{lemma}
\begin{remark}
The significance of $\Psi$ is explained in Remark~\ref{psirem}.
\end{remark}

%\begin{remark}
%If we set $\Gamma = \Sigma_w-\delta I$, for $\delta >0$, then $P(A_k) \to 1$ using the fact that increments $w_t$ are independent and the law of large numbers.
%\end{remark}

\begin{proof}
Let $\e>0$ and Observe that
\begin{align*}
 \sup_{\theta' \in B(\theta,\e)} R_T^\pi(\theta')\geq\E_{\theta'\sim \lambda} R_T^\pi(\theta')
\end{align*}
for any prior $\lambda \in C^\infty_c \left(B(\theta,\e)\right)$. First, we reduce this problem to a sequence of Bayesian estimation problems on a variable scale $\e$\footnote{This scale must be decreasing, but not too fast due certain constraints imposed by the arguments used in Lemma~\ref{informationcomparisonlemma}.}. Second, we decouple the resulting information matrix from the underlying state process by a simple truncation argument.

\paragraph{Reduction to Bayesian estimation.} Expanding the summand (\ref{ziemannregexp}) of Theorem~\ref{regexpthm}
\begin{multline*}
 \E_{\theta'\sim \lambda}^\pi\left( (u_t-K(\theta')x_t)^\top\big[R +B^\top(\theta') P(\theta') B(\theta')\big](u_t-K(\theta')x_t)\right)\\
 =\E_{\theta'\sim \lambda}^\pi\tr\left(\big[R +B^\top(\theta') P(\theta') B(\theta')\big](u_t-K(\theta')x_t) (u_t-K(\theta')x_t)^\top\right)\\
  \geq\tr \left(\big[R +B^\top(\theta) P(\theta) B( \theta)+O(\e)\big]\E^\pi_{\theta'\sim \lambda} \left[(u_t-K(\theta')x_t) (u_t-K(\theta')x_t)^\top \right]\right)
  \end{multline*}
  by continuity.  It follows that
 \begin{multline}
 \label{td:1}
 \sup_{\theta' \in B(\theta,\e T^{-1/4})} R_T^\pi(\theta') \geq \sum_{k=1}^{\lfloor\sqrt{T-1}\rfloor} \sum_{t=k\lceil\sqrt{T}\rceil}^{(k+1)\lfloor\sqrt{T}\rfloor}\tr \Bigg(\big[R +B^\top(\theta) P(\theta) B(\theta)+O(\e)\big]\\ \times\E^\pi_{\theta'\sim \lambda} \left[(u_t-K(\theta')x_t) (u_t-K(\theta')x_t)^\top \right]\Bigg)  +O(1)\\
=
\sum_{k=1}^{\lfloor\sqrt{T-1}\rfloor} \sum_{t=k\lceil\sqrt{T}\rceil}^{(k+1)\lfloor\sqrt{T}\rfloor}\tr \Bigg(\big[R +B^\top(\theta) P(\theta) B(\theta)+O(\e)\big]\\ \times\E^\pi_{\theta'\sim \lambda} \E^\pi_{\theta'\sim \lambda} \left[(u_t-K(\theta')x_t) (u_t-K(\theta')x_t)^\top \Big| \mathcal{S}_k \right]\Bigg)+O(1)
 \end{multline} 
 where we choose $\mathcal{S}_{k,T}$ to be the $\sigma$-field generated by $(x_{k\lceil\sqrt{T}},\dots, x_{(k+1)\lceil\sqrt{T}})$. We now apply Van Trees' inequality (\ref{VTineq}) to the regular conditional density in the interior expectation in (\ref{td:1}). Denote by $\I(\theta|\mathcal{S}_{k,T})$ the Fisher information corresponding to the density $p(x_1,\dots,x_{k\lceil\sqrt{T}-1}| \theta,\mathcal{S}_{k,T})$ and $\J(\lambda | \mathcal{S}_{k,T})$ the location information corresponding to the density of $\theta$ given $\mathcal{S}_k$. Then by (\ref{VTineq}) we have
 \begin{multline}
 \label{crusedhere}
  \E^\pi_{\theta'\sim \lambda} \left[(u_t-K(\theta')x_t) (u_t-K(\theta')x_t)^\top \Big| \mathcal{S}_{k,T} \right]\\
  \succeq \left(\E^\pi_{\theta'\sim \lambda} [\dop_\theta K(\theta')x_t |\mathcal{S}_{k,T} ] \right)^\top\left( \E^\pi_{\theta'\sim \lambda} \I (\theta|\mathcal{S}_{k,T})|\mathcal{S}_{k,T}] + \J(\lambda | \mathcal{S}_{k,T})\right)^{-1}\left(\E^\pi_{\theta'\sim \lambda} [\dop_\theta K(\theta')x_t |\mathcal{S}_{k,T} ] \right).
 \end{multline}

 \paragraph{Decoupling by truncation.} The issue with (\ref{crusedhere}) is that the state process $x_t$ is correlated with the information term, due to the conditioning. To remedy this, we combine (\ref{crusedhere}) with (\ref{td:1}) using Lemma~\ref{kronprodlemma} and introduce the truncation events $A^T_k = \{\sum_{t=k\lceil\sqrt{T}\rceil}^{(k+1)\lfloor\sqrt{T}\rfloor}x_t x_t^\top \succeq \sqrt{T} \Gamma \}$ for $\Gamma\succ 0$. We find
 \begin{multline*}
 \sup_{\theta' \in B(\theta,\e T^{-1/4})} R_T^\pi(\theta') \geq O(1) +  \sum_{k=1}^{\lfloor\sqrt{T-1}\rfloor} \sum_{t=k\lceil\sqrt{T}\rceil}^{(k+1)\lfloor\sqrt{T}\rfloor}\tr \Bigg(\left(x_t x_t^\top \otimes \big[R +B^\top(\theta) P(\theta) B(\theta)+O(\e)\big]\right)\\ \times\E^\pi_{\theta'\sim \lambda} \left( \E^\pi_{\theta'\sim \lambda}[\dop_\theta \VEC K(\theta)|\mathcal{S}_{k,T}] \right) \left(\left( \E^\pi_{\theta'\sim \lambda} \I (\theta|\mathcal{S}_{k,T})|\mathcal{S}_{k,T}] + \J(\lambda | \mathcal{S}_{k,T})\right)^{-1}  \right) \left( \E^\pi_{\theta'\sim \lambda}[\dop_\theta \VEC K(\theta)|\mathcal{S}_{k,T}] \right)^\top\Bigg)\\
 % inf over K(\theta)
 \geq
  O(1) +  \sum_{k=1}^{\lfloor\sqrt{T-1}\rfloor} \sum_{t=k\lceil\sqrt{T}\rceil}^{(k+1)\lfloor\sqrt{T}\rfloor}\tr \Bigg(\left(x_t x_t^\top \otimes\big[R +B^\top(\theta) P(\theta) B(\theta)+O(\e)\big]\right)\\ 
  \times \left(\dop_\theta \VEC K(\theta)+O(\e) \right)\E^\pi_{\theta'\sim \lambda}\left( \left( \E^\pi_{\theta'\sim \lambda} \I (\theta|\mathcal{S}_{k,T})|\mathcal{S}_{k,T}] + \J(\lambda | \mathcal{S}_{k,T})\right)^{-1}  \right) \left(\dop_\theta \VEC K(\theta)+O(\e)\right)^\top \Bigg)\\
  % truncate and switch summation
   \geq
  O(1) +  \sum_{k=1}^{\lfloor\sqrt{T-1}\rfloor}\tr \Bigg(\left( \sum_{t=k\lceil\sqrt{T}\rceil}^{(k+1)\lfloor\sqrt{T}\rfloor}  x_t x_t^\top \otimes\big[R +B^\top(\theta) P(\theta) B(\theta)+O(\e)\big]\right)
\left(\dop_\theta \VEC K(\theta)+O(\e) \right) \\
\times  \E^\pi_{\theta'\sim \lambda}\mathbf{1}_{A_k}\left(\left( \E^\pi_{\theta'\sim \lambda} \I (\theta|\mathcal{S}_{k,T})|\mathcal{S}_{k,T}] + \J(\lambda | \mathcal{S}_{k,T})\right)^{-1}  \right)\left(\dop_\theta \VEC K(\theta)+O(\e)\right)^\top \Bigg)\\
% execute truncation
 \geq
  O(1) +  \sum_{k=1}^{\lfloor\sqrt{T-1}\rfloor}\tr \Bigg(\left(\Gamma \otimes  \big[R +B^\top(\theta) P(\theta) B(\theta)+O(\e)\big]\right)
\left(\dop_\theta \VEC K(\theta)+O(\e) \right)  \\
\times  \E^\pi_{\theta'\sim \lambda}\left( \sqrt{T}  \mathbf{1}_{A_k} \left( \E^\pi_{\theta'\sim \lambda} \I (\theta|\mathcal{S}_{k,T})|\mathcal{S}_{k,T}] + \J(\lambda | \mathcal{S}_{k,T})\right)^{-1}  \right)\left(\dop_\theta \VEC K(\theta)+O(\e)\right)^\top\Bigg)\\
% Jensen
 \geq
  O(1) +  \sum_{k=1}^{\lfloor\sqrt{T-1}\rfloor}\tr \Bigg(\left( \Gamma \otimes \big[R +B^\top(\theta) P(\theta) B(\theta)+O(\e)\big]\right)
\left(\dop_\theta \VEC K(\theta)+O(\e) \right)\\
\times  [\mathbf{P}^\pi_{\theta'\sim \lambda}( A_k)]^2 \left( \sqrt{T}  \left( \E^\pi_{\theta'\sim \lambda}\left[\mathbf{1}_{A_k}( \E^\pi_{\theta'\sim \lambda} \I (\theta|\mathcal{S}_{k,T})|\mathcal{S}_{k,T}] + \J(\lambda | \mathcal{S}_{k,T}) )\right]\right)^{-1}  \right)\left(\dop_\theta \VEC K(\theta)+O(\e)\right)^\top  \Bigg)\\
 \geq
  O(1) +  \sum_{k=1}^{\lfloor\sqrt{T-1}\rfloor} [\mathbf{P}^\pi_{\theta'\sim \lambda}( A_k)]^2\tr \Bigg(\left(\Gamma \otimes \big[R +B^\top(\theta) P(\theta) B(\theta)+O(\e)\big]\right)
\left(\dop_\theta \VEC K(\theta)+O(\e) \right) \\
\times   \left( \sqrt{T} 
 \left( \E^\pi_{\theta'\sim \lambda} \I^{(k+1)\lfloor\sqrt{T}\rfloor} (\theta) + \J(\lambda )\right)^{-1}  \right)\left(\dop_\theta \VEC K(\theta)+O(\e)\right)^\top \Bigg)
 \end{multline*}
 where the penultimate step follows by Matrix Jensen's inequality and the final step by the chain rule for Fisher information. The result follows by re-arranging and adding the positive definite pertubation $\Psi\succeq 0$, (which technically makes the inequality weaker).
\end{proof}

\begin{remark}
\label{gammasigmaremark}
For the choice of $\Gamma = \Sigma_w$ in the Lemma~\ref{decouplemma} we have $\lim_T P(A^T_k) \to 1$. To see this note that for any length $\tau$ trajectory
\begin{align*}
\frac{1}{\tau}\sum_{t=\tau_0}^{\tau + \tau_0} x_t x_t^\top &= \frac{1}{\tau}\sum_{t=\tau_0}^{\tau + \tau_0} (Ax_{t-1}+Bu_{t-1}+w_{t-1})(Ax_{t-1}+Bu_{t-1}+w_{t-1})^\top\\
&=\frac{1}{\tau}\sum_{t=\tau_0}^{\tau + \tau_0} (Ax_{t-1}+Bu_{t-1})(Ax_{t-1}+Bu_{t-1})^\top +\frac{1}{\tau}\sum_{t=\tau_0}^{\tau + \tau_0} w_{t-1} w_{t-1}^\top \\
&+\frac{1}{\tau}\sum_{t=\tau_0}^{\tau + \tau_0} (Ax_{t-1}+Bu_{t-1}) w_{t-1}^\top + \frac{1}{\tau}\sum_{t=\tau_0}^{\tau + \tau_0}w_{t-1} (Ax_{t-1}+Bu_{t-1})^\top\\
&\succeq \frac{1}{\tau}\sum_{t=\tau_0}^{\tau + \tau_0} w_{t-1} w_{t-1}^\top+\frac{1}{\tau}\sum_{t=\tau_0}^{\tau + \tau_0} (Ax_{t-1}+Bu_{t-1}) w_{t-1}^\top + \frac{1}{\tau}\sum_{t=\tau_0}^{\tau + \tau_0}w_{t-1} (Ax_{t-1}+Bu_{t-1})^\top.
\end{align*}
The first term on the RHS above converges almost surely to $\Sigma_w$ by the law of large numbers, whereas the second and third terms converge to zero by the martingale law of large numbers. This establishes the claim.
\end{remark}

The following lemma details the pertubation theory of Fisher information for low regret policies.

\begin{lemma}
\label{informationcomparisonlemma}
%Assume that 
%\begin{align*}
%\limsup_{n\to \infty} \sup_{\theta' \in B(\theta,\e)} \frac{R_n^\pi(x,\theta')}{\sqrt{n}} \leq C_\pi
%\end{align*}
%for some $C_\pi>0$ and $\e >0 $. 
Assume A1-A4. Fix $\e>0$ and a of prior density $\lambda$ supported on $B(\theta, \e)$. Fix a positive semidefinite matrice, $\Psi\succeq 0 $ and let
\begin{align*}
\E^\pi_{\theta\sim \lambda}\left[ \I^{t}(\theta)\right]+\J(\lambda) +\Psi=  V_1 \Lambda_1 V_1^\top + V_0 \Lambda_0 V_0^\top
\end{align*}
where $\Lambda_1= \diag(\lambda_1,\dots, \lambda_{d_\theta - \dim \mathtt{U}})$, $\Lambda_0=\diag(\lambda_{d_\theta - \dim \mathtt{U}+1},\dots, \lambda_{d_\theta})$ and $V_1, V_0$ spanning the corresponding eigenspaces of dimension $d_\theta - \dim \mathtt{U} $ and  $\dim \mathtt{U}$ respectively. Let further
\begin{align*}
\E^\pi_{\theta\sim \lambda}\sum_{j=0}^{t}  L^\top (K(\theta)x_j x_j K^\top(\theta)  \otimes \J(p)  ) L= \tilde V_1 \tilde \Lambda_1\tilde V_1^\top + \tilde V_0 \tilde \Lambda_0  \tilde V_0^\top
\end{align*} 
with the variables $\tilde V, \tilde \Lambda$ defined similarly. If further $\Psi\tilde V_0 =0$ and $ \sigma_{\min} (\tilde V_1^\top \Psi \tilde V_1 )>0$, then the small part of the spectrum $\Lambda_0$ satisfies the inequality
\begin{multline}
\label{lambdabound}
\|\Lambda_0\|_\infty \leq \|L\|_\infty^2  \|Q^{-1}\|_\infty \left(V_{t}^{\pi^*}(\theta)+R_{t}^{\pi}(\theta)\right)  \tr [J(p)] \|\dop^2_\theta \tr K(\theta')K^\top(\theta') \|_\infty\times \frac{\e^2}{2}\\
+2 \| L\|^2_\infty \| [R+B^\top(\theta) K(\theta) B(\theta)]^{-1}\|_\infty \tr \left( \J(p)  \right)  R_{t}^\pi(\theta)\\
 + \|\J(\lambda)\|_\infty+O\left( \e^3 \left[V_{t-\tau}^{\pi^*}(\theta)+R_{t}^{\pi}(\theta)\right]\right)=:\Delta_{t,T}
\end{multline}
and the $k=1,\dots,\dim \mathtt{U}$ smallest singular values in particular satisfy
\begin{multline}
\sigma_{k}( \Lambda_0) \leq  \|L\|_\infty^2  \|Q^{-1}\|_\infty \left(V_{t}^{\pi^*}(\theta)+R_{t}^{\pi}(\theta)\right) \frac{ \tr \J(p)}{k} \|\dop^2_\theta \tr K(\theta')K^\top(\theta') \|_\infty\times \frac{\e^2}{2}\\
+2 \| L\|^2_\infty \| [R+B^\top(\theta) K(\theta) B(\theta)]^{-1}\|_\infty  \frac{\tr\J(p)}{k} R_{t}^\pi(\theta)\\
 + \|\J(\lambda_t)\|_\infty+O\left( \e^3 \left[V_{t}^{\pi^*}(\theta)+R_{t}^{\pi}(\theta)\right]\right).
 \end{multline}
Finally, the eigenspaces corresponding to the small eigenvalues satisfy the following perturbation inequality
\begin{align}
d_\infty(V_0,\tilde V_0) \leq \frac{\Delta_{t,T}}{ \sigma_{\min} (\tilde V_1^\top \Psi \tilde V_1 )}
\end{align}
where the metric $d_\infty$ is defined in (\ref{ddef}).
\end{lemma}

\begin{remark}
\label{psirem}
$\Psi$ is essentially a dummy variable, which we will send to infinity on the subspace orthogonal to $\tilde V_0$ ($\tilde V_0$ spans $\mathtt U$). Intuitively, this concentrates the prior on the  information singular subspace.
\end{remark}

\begin{proof}
We are interested in the singular perturbation theory of $\E^\pi_{\theta\sim \lambda}\left[ \I^{t}(\theta)\right]+\J(\lambda)$ under the hypothesis that regret is $O(\sqrt{T})$. Our main tool will be the Davis-Kahan $\sin\theta$ Theorem, for which we now set the stage.  Observe  that for any matrix $K$
\begin{multline}
\label{kxexp}
u_j u_j^\top = (Kx_j +u_j-Kx_j)(Kx_j+u_j-Kx_j)^\top \\
= Kx_jx_j^\top K +(u_j-Kx_j)(u_j-Kx_j)^\top 
+ Kx_j(u_j^\top-x_j^\top K^\top)+(u_j-Kx_j)x_j^\top\\
\preceq 2 Kx_jx_j^\top K+2(u_j-Kx_j)(u_j-Kx_j)^\top.
\end{multline}
In particular we may write, pointwise for each $\theta$,
\begin{multline}
\label{fisherkxexp}
 \I^{t}(\theta) =\E \sum_{j=0}^{t}  L^\top (u_ju_j^\top  \otimes \J(p)  ) L\\
 =\E \sum_{j=0}^{t}  L^\top (Kx_j x_j K^\top  \otimes \J(p)  ) L+\E \sum_{j=0}^{t}  L^\top ((u_j-Kx_j)(u_j-Kx_j)^\top  \otimes \J(p)  ) L\\
 +\E \sum_{j=0}^{t}  L^\top (Kx_j(u_j-Kx_j)^\top  \otimes \J(p)  ) L+\E \sum_{j=0}^{t}  L^\top ((u_j-Kx_j)x_j^\top K^\top  \otimes \J(p)  ) L\\
 \preceq 2\left(\E \sum_{j=0}^{t}  L^\top (Kx_j x_j K^\top  \otimes \J(p)  ) L+\E \sum_{j=0}^{t}  L^\top ((u_j-Kx_j)(u_j-Kx_j)^\top  \otimes \J(p)  ) L\right)
\end{multline}
where the equality in (\ref{fisherkxexp}) follows from the equality in (\ref{kxexp}) and the inequality from the corresponding inequality and the fact that Kronecker products preserve the semidefinite order.

For the duration of this proof, it will be convenient to denote
\begin{align*}
N &= \E^\pi_{\theta\sim \lambda}\I^{t}(\theta)+\rho\tau M+\J(\lambda)-\Psi ,\\
\tilde N&=\E^\pi\sum_{j=0}^{t}  L^\top (K(\theta)x_j x_j^\top K^\top(\theta)  \otimes \J(p)  ) L ,\\
T' &= \E^\pi\sum_{j=0}^{t}  L^\top (K(\theta)x_j x_j^\top K^\top(\theta)  \otimes \J(p)  ) L-\rho\tau M-\J(\lambda) -\Psi\\
 &-\E^\pi_{\theta \sim \lambda}\sum_{j=0}^{t}  L^\top (K(\theta)x_j x_j^\top K^\top(\theta)  \otimes \J(p)  ) L\\
 &-\E^\pi_{\theta\sim \lambda} \sum_{j=0}^{t}  L^\top ((u_j-K(\theta)x_j)(u_j-K(\theta)x_j)^\top  \otimes \J(p)  ) L\\
 &-\E^\pi_{\theta\sim \lambda} \sum_{j=0}^{t}  L^\top (K(\theta)x_j(u_j-K(\theta)x_j)^\top  \otimes \J(p)  ) L\\
 &-\E^\pi_{\theta\sim \lambda} \sum_{j=0}^{t}  L^\top ((u_j-K(\theta)x_j)x_j^\top K^\top (\theta) \otimes \J(p)  ) L
\end{align*}
In this way, $N$ is a pertubation of $\tilde N$ of magnitude $T'$: $\tilde N = N +T'$ and we think of $\tilde N$ as Fisher information corresponding to the optimal\footnote{Although, note that the expectation is with respect to $\pi$, not $\pi^*$.} linear feedback $K$ at the nominal instance $\theta$. The rest of the proof consists of showing that $T'$ cannot be ``too large'' on the nullspace of $\tilde N$. Recall also the spectral decompositions
\begin{align}
\begin{cases}
\label{specdecs}
N = V_1 \Lambda_1 V_1^\top + V_0 \Lambda_0 V_0^\top\\
\tilde N =  \tilde V_1 \tilde \Lambda_1\tilde V_1^\top + \tilde V_0 \tilde \Lambda_0  \tilde V_0^\top
\end{cases}
\end{align}
where $\Lambda_1= \diag(\lambda_1,\dots, \lambda_{g})$, $\Lambda_0=\diag(\lambda_{g+1},\dots, \lambda_{d_\theta}), g=d_\theta - \dim \mathtt{U}$ and $V_1, V_0$ partial isometries and similarly for the variables in $\tilde V_i,\tilde \Lambda_i$. It is clear by construction of these matrices that they are positive semidefinite so that $\lambda_i \geq 0, \tilde \lambda_i \geq 0$. Further, we may take $V_0, \tilde V_0$ to have orthonormal columns and span the eigenspaces corresponding to the smallest $d_\theta-g$ eigenvalues of $N, \tilde N$ respectively. Take now $g = \dim \ker \tilde N$ in which case we observe that $\tilde V_0$ spans $\ker \tilde N$. To apply the Davis-Kahan $\sin\theta$-Theorem, we need to prove that $\|T' \tilde V_0\|$ is small and that there exist constants $\alpha\geq 0, \delta > 0$ such that the seperation conditions $\sigma_{\min}(\Lambda_1) \geq \alpha + \delta$ and $\sigma_{\max}(\tilde \Lambda_0) \leq \alpha$ hold. Note now that $\tilde \Lambda_0=0$ so that we may take $\alpha=0$ and $\delta= \sigma_{\min} \tilde V_1^\top \Psi \tilde V_1$. Thus we are left with proving a bound on $\|T' \tilde V_0\|$.

We now procede with a bound on $\|T' \tilde V_0\|$.  Observe by (\ref{fisherkxexp}) that
\begin{multline}
0\succeq -\E^\pi_{\theta' \sim \lambda}\sum_{j=0}^{t}  L^\top (K(\theta')x_j x_j K^\top(\theta')  \otimes \J(p)  ) L \\
-\E^\pi_{\theta'\sim \lambda} \sum_{j=0}^{t }  L^\top ((u_j-K(\theta')x_j)(u_j-K(\theta')x_j)^\top  \otimes \J(p)  ) L\\
 -\E^\pi_{\theta'\sim \lambda} \sum_{j=0}^{t }  L^\top (K(\theta')x_j(u_j-H(\theta')x_j)^\top  \otimes \J(p)  ) L-\E^\pi_{\theta'\sim \lambda_t} \sum_{j=0}^{t }  L^\top ((u_j-H(\theta')x_j)x_j^\top H^\top (\theta') \otimes \J(p)  ) L\\
 \succeq -2\underbrace{\E^\pi_{\theta'\sim \lambda} \sum_{j=0}^{t }  L^\top (K(\theta')x_j x_j K^\top(\theta')  \otimes \J(p)  ) L}_{=:T_1'}\\
 -2\underbrace{\E^\pi_{\theta'\sim \lambda_t} \sum_{j=0}^{t }  L^\top ((u_j-K(\theta')x_j)(u_j-K(\theta')x_j)^\top  \otimes \J(p)  ) L}_{=:T_2'}.
\end{multline}
We now prove separetely that $T_1'$  and $T_2'$ are small on $\ker \tilde N$. It is helpful to note that $T_1$ is proportional to the $\lambda$-integral of $\tilde N=\tilde N(\theta)$. Together, this will imply that the nullspace of $\tilde N$ is approximately that of $N$, provided that $\rho\tau M+\J(\lambda)$ is not too large.

We first prove that $T_1'\tilde V_0$ is small. To do so, note that since $\tilde N \tilde V_0 =0$ we may write, by Taylor expansion,
\begin{multline}
\label{t1v0taylorexp}
 \tr \E^\pi\left[ \sum_{j=0}^{t } L^\top (K(\theta')x_j x_j K^\top(\theta')  \otimes \J(p)  ) L\tilde V_0   \Big| \theta'\right]\\ = \frac{1}{2} (\theta-\theta')^\top \mathcal{H}_{t } (\theta-\theta') + \textnormal{ Higher order terms}
\end{multline}
where
\begin{align*}
\mathcal{H}_{t } = \dop_\theta^2 \tr \E^\pi\left[ \sum_{j=0}^{t } L^\top (K(\theta')x_j x_j K^\top(\theta')  \otimes \J(p)  ) L\tilde V_0\Big|\theta'\right].
\end{align*}
In the Taylor expansion above, the zeroth and first order terms are zero since $\theta$ is a local minimum due to $\tilde N \tilde V_0 =0$ and that $\tilde T_1$ is positive semidefinite. Now the support of $\lambda$ is contained in $B(\theta, \e)$,  wherefore we in fact have that
\begin{align*}
 \tr T_1'\tilde V_0  \leq \|\mathcal{H}_{t }\|_\infty \times \frac{\e^2}{2} +O\left(\|\dop_\theta \mathcal{H}_{t }\|_\infty\e^3\right)
\end{align*}
which is the main quantity to bound, and which drives our result, in combination with the observation that
\begin{align*}
\| T_1'\tilde V_0\|_\infty \leq  \tr T_1'\tilde V_0   \textnormal{ and }\sigma_{k}(T_1'\tilde V_0) \leq \frac{1}{k} \tr T_1'\tilde V_0 .
\end{align*}

Next we prove that the derivatives of (\ref{t1v0taylorexp}) grow at most linearly.\footnote{Doing so later enables us to conclude that $\| T_1'\tilde V_0\|_\infty = O(\sqrt{t})$.} To do so, note that since we are conditioning on $\theta'$, we have by linearity
\begin{multline*}
 \tr \E^\pi\left[ \sum_{j=0}^{t } L^\top (K(\theta')x_j x_j K^\top(\theta')  \otimes \J(p)  ) L\tilde V_0\Big|\theta'\right]\\
=\tr  L^\top \left(K(\theta')\E^\pi\left[ \sum_{j=0}^{t }  x_j x_j\right] K^\top(\theta')  \otimes \J(p) \right )  L\tilde V_0
\end{multline*}
and so, by smoothness of $K(\cdot)$, it suffices to bound the term within the expectation above. Note now that 
\begin{multline*}
\E^\pi\left[ \sum_{j=0}^{t }  x_j x_j\right] \leq \|Q^{-1}\|_\infty \E^\pi\left[ \sum_{j=0}^{t }  x_j Q x_j\right] \leq \|Q^{-1}\|_\infty V_{t }^\pi(\theta)\\
=  \|Q^{-1}\|_\infty \left(V_{t }^{\pi^*}(\theta)+R_{t }^{\pi}(\theta)\right)
\end{multline*}
using the positive definitess of $Q$ and we thus have that the derivatives of (\ref{t1v0taylorexp}) are $O\left(V_{t }^{\pi^*}(\theta)+R_{t }^{\pi}(\theta)\right)$. In fact, we have the estimate
\begin{align*}
\|\mathcal{H}_{t }\|_\infty \leq  2 \|L\|_\infty^2  \|Q^{-1}\|_\infty \left(V_{t }^{\pi^*}(\theta)+R_{t }^{\pi}(\theta)\right)  \tr [J(p)] \|\dop^2_\theta \tr K(\theta')K^\top(\theta') \|_\infty
\end{align*}
so that 
\begin{multline}
\label{t1v1finalbound}
\| T_1'\tilde V_0\|_\infty \leq  \|L\|_\infty^2  \|Q^{-1}\|_\infty \left(V_{t }^{\pi^*}(\theta)+R_{t }^{\pi}(\theta)\right)  \tr [J(p)] \|\dop^2_\theta \tr K(\theta')K^\top(\theta') \|_\infty\times \frac{\e^2}{2}\\ +O\left( \e^3 \left[V_{t }^{\pi^*}(\theta)+R_{t }^{\pi}(\theta)\right]\right).
\end{multline}

Next, we prove that $T_2'$ is small. In fact, this is an immediate consequence of having low regret. We have to estimate
\begin{align*}
2\left\|\E \sum_{j=0}^{t }  L^\top ((u_j-Kx_j)(u_j-Kx_j)^\top  \otimes \J(p)  ) L\right\|_\infty.
\end{align*}
It suffices to estimate
\begin{multline}
\begin{aligned}
\label{selfbounding}
\tr T_2'\tilde V_0 &\leq 2 \| L\|^2_\infty \sum_{j=0}^{t }\tr\left(\E ((u_j-Kx_j)(u_j-Kx_j)^\top  \otimes \J(p)  )\right)\\
&= 2 \| L\|^2_\infty \sum_{j=0}^{t }\tr\left(\E ((u_j-Kx_j)(u_j-Kx_j)^\top\right)\tr \left( \J(p)  \right)\\
&= 2 \| L\|^2_\infty \tr  \J(p) \\
&\times \sum_{j=0}^{t }\tr\left(\E ((u_j-Kx_j)[R+B^\top(\theta) K(\theta) B(\theta)][R+B^\top(\theta) K(\theta) B(\theta)]^{-1}(u_j-Kx_j)^\top\right)\\
&\leq 2 \| L\|^2_\infty \| [R+B^\top(\theta) K(\theta) B(\theta)]^{-1}\|_\infty \tr \left( \J(p)  \right)  R_{t }^\pi(\theta).
\end{aligned}
\end{multline}
We thus conclude that
\begin{multline}
\begin{aligned}
\label{Tapprox}
\| T\tilde V_0\|_\infty &\leq \|T'_1\tilde V_0\|_\infty+\|T'_2\tilde V_0\|_\infty +\rho\tau\| M\|_\infty+ \|\J(\lambda)\|_\infty\\
&\leq  \|L\|_\infty^2  \|Q^{-1}\|_\infty \left(V_{t }^{\pi^*}(\theta)+R_{t }^{\pi}(\theta)\right)  \tr [J(p)] \|\dop^2_\theta \tr K(\theta')K^\top(\theta') \|_\infty\times \frac{\e^2}{2}\\
&\:\:\:\:+2 \| L\|^2_\infty \| [R+B^\top(\theta) K(\theta) B(\theta)]^{-1}\|_\infty \tr \left( \J(p)  \right)  R_{t }^\pi(\theta)\\
&\:\:\:\: + \|\J(\lambda_t)\|_\infty+O\left( \e^3 \left[V_{t }^{\pi^*}(\theta)+R_{t }^{\pi}(\theta)\right]\right)
 \end{aligned}
\end{multline}
and also obtain the following estimate for the smallest singular value
\begin{multline}
\begin{aligned}
\label{Tapprox2}
\sigma_{k}( T \tilde V_0)& \leq (k)^{-1} \left[\tr(T'_1\tilde V_0)+\tr(T'_2\tilde V_0)\right] +\rho\tau\| M\|_\infty+ \|\J(\lambda)\|_\infty\\
&\leq  \|L\|_\infty^2  \|Q^{-1}\|_\infty \left(V_{t }^{\pi^*}(\theta)+R_{t }^{\pi}(\theta)\right)  \tr \frac{J(p)}{k} \|\dop^2_\theta \tr K(\theta')K^\top(\theta') \|_\infty\times \frac{\e^2}{2}\\
&\:\:\:\:+2 \| L\|^2_\infty \| [R+B^\top(\theta) K(\theta) B(\theta)]^{-1}\|_\infty \tr \frac{J(p)}{k} R_{t }^\pi(\theta)\\
 &\:\:\:\:+ \|\J(\lambda_t)\|_\infty+O\left( \e^3 \left[V_{t }^{\pi^*}(\theta)+R_{t }^{\pi}(\theta)\right]\right)
 \end{aligned}
\end{multline}
using (\ref{t1v1finalbound}), (\ref{selfbounding}), the fact that $\tilde V_0$ spans $\ker \tilde N$, $\|\tilde V_0\|_\infty \leq 1$, and a triangle inequality.\end{proof}

\section{Proof of Theorem~\ref{regexpthm}: Regret, Riccati and Bellman}
\label{bkappendix}
Consider the long-run average cost
\begin{align}
\label{ac}
g^\pi(\theta) =\limsup_{T\to\infty} \frac{1}{T} \sum_{t=1}^T \E^\pi_{\theta} c(x_t,u_t).
\end{align}
It can be shown that there exists a policy $\pi^{\infty}(\theta)$ minimizing $(\ref{ac})$ that solves the following Bellman equation for each $x$
\begin{align}
\label{acoe}
\lambda(\theta) + h(x,\theta) = \min_{u\in \mathbb{R}^{d_u}} \{ c(x,u)+ \E_w[h(A(\theta)x+B(\theta)u+w,\theta)]\}
\end{align}
where $c(x,u) = x^\top Q x+ u^\top R u $ as before, $\lambda(\theta) = \inf_\pi g^\pi(\theta)=\tr (P(\theta)\Sigma_w)$ and $h(x,\theta) = x^\top P(\theta) x$ where $P(\theta)$ solves the following Riccati equation
\begin{align}
\label{ARE}
P (\theta)  = A^\top(\theta)[P(\theta) - P(\theta) B(\theta)(R+B^\top(\theta)P(\theta)B(\theta))^{-1}B^\top(\theta)P(\theta)]A(\theta) + Q.
\end{align}
Our analysis rests on the following quantity, which measures the suboptimality of any policy with respect to the gap in the Bellman equation (\ref{acoe})
\begin{multline}
\label{phifunc}
\phi(x,u,\theta)  =   c(x,u)+ \E_w[h(A(\theta)x+B(\theta)u+w,\theta)] - \lambda(\theta) - h(x,\theta) \\ 
=x^\top Q x + u^\top R u + \E_w [(A(\theta)x+B(\theta)u+w)(P(\theta))(A(\theta)x+B(\theta)u+w)]-\lambda(\theta)-x^\top P(\theta) x\\
=x^\top Q x + u^\top R u +  [A(\theta)x+B(\theta)u]P(\theta)[A(\theta)x+B(\theta)u)]-x^\top P(\theta) x
\end{multline}
Observe that (\ref{acoe}) may be written using (\ref{phifunc}) as $\min_{u\in \mathbb{R}^{d_u}} \phi(x,u,\theta)=0$. Our goal will be to prove that regret is approximately a cumulative sum over the Bellman sub-optimalities $\phi(x,u,\theta)\neq 0$ induced by policy $\pi \neq \pi^*$.

We work on the stability region $\{\theta' \in \mathbb{R}^{d_\theta} |\rho(A-B(\theta')K(\theta)) \leq 1 -\zeta\}$, $\zeta \in (0,1)$. The following lemma, originally proven by \cite{burnetas1997optimal} in the context of finite state and action space Markov Decision Processes, gives an expression for regret in terms of the function $\phi$ in (\ref{phifunc}). 
\begin{lemma}
\label{testqlemma}
Under assumptions A1-A4, for any policy $\pi$, it holds that
\begin{align}
\label{bkregexp}
R_T^\pi(\theta) = \sum_{t=0}^{T-1} \E_{\theta}^\pi \phi(x_t,u_t,\theta)+O(1)
\end{align}
where the function $\phi$ is given in (\ref{phifunc}) and $h$ is as in (\ref{acoe}). The term $O(1)$ can be kept uniformly bounded on the stability region $\{\theta' \in \mathbb{R}^{d_\theta} |\rho(A-B(\theta')K(\theta)) \leq 1 -\zeta\}$, $\zeta \in (0,1)$.
\end{lemma}

\begin{proof}
By (\ref{asymptoticequivalence}), we may write
\begin{align}
\label{riccatiregexp}
\sum_{t=0}^{T-1}\E_{x,\theta}^{\pi^*} c(x_t,u_t) = Tg(\theta) +h(x,\theta) +O(1).
\end{align}
Hence $Tg(\theta) +h(x,\theta) = \sum_{t=0}^{T-1}\E_{x,\theta}^{\pi^*} c(x_t,u_t)+O(1) $. That is, $V_T^*(x,\theta) = Tg(\theta)+h(x,\theta)+O(1)$.  Substituting this into the definition of regret gives $R_T^\pi(x,\theta) = V_T^\pi(x,\theta) -Tg(\theta)-h(x,\theta)+O(1)$. For convenience set $D_T^\pi(x,\theta) =V_T^\pi(x,\theta) -T\lambda(\theta)-h(x,\theta) $.  Now, express $V_T^\pi(x,\theta)$ recursively as $V_T^\pi(x,\theta) = \E_{x,\theta}^\pi [\E_{u_0}[\E_{x_1}c(x,u_0)+V_{T-1}^\pi(x_1,\theta)]]].$ Combining  yields
\begin{multline}
D_T^\pi(x,\theta) = \E_{x,\theta}^\pi c(x,u_0)+V_{T-1}^\pi(x_1,\theta)-Tg(\theta)-h(x,\theta)\\
=\E_{x,\theta}^\pi c(x,u_0)+\E^u_{x,\theta} h(x_1,\theta)-g(\theta)-h(x,\theta)+V_{T-1}^\pi(x_1,\theta)-(T-1)g(\theta)-\E^u_{x,\theta} h(x_1,\theta)\\
=\E_{x,\theta}^\pi \phi(x,u,\theta)+D_{T-1}^\pi(x,\theta)=\dots=\sum_{t=0}^{T-1} \E_{x,\theta}^\pi \phi(x_t,u_t,\theta).
\end{multline}
The result follows since $R_n^\pi(x,\theta) = D_n^\pi(x,\theta)+\E_{x,\theta}^{\pi^*}h(x_n,\theta)$.
\end{proof}

We 

\paragraph{Proof of Theorem~\ref{regexpthm}}
 We have that 
\begin{align*}
\dop_u^2 \phi(x,u,\theta) = 2R +2B^\top(\theta) P(\theta) B(\theta).
\end{align*}
We now taylor expand $\phi(x,u,\theta)$ at $u= K(\theta)x$. Since this is a quadratic function with global minimum $\phi(x,K(\theta)x,\theta)=0$ by optimality of $K(\theta)x$, it follows that
\begin{align}
\label{phitaylorexp}
\phi(x,u,\theta) =  (u-K(\theta)x)^\top\big[R +B^\top(\theta) P(\theta) B(\theta)\big](u-K(\theta)x).
\end{align}
Substituting (\ref{phitaylorexp}) into (\ref{bkregexp}) we obtain (\ref{ziemannregexp}). $\hfill \blacksquare$

The proof of Lemma~\ref{testqlemma} essentially follows verbatim from  \cite{burnetas1997optimal} and once this lemma is established Theorem~\ref{regexpthm} is an easy consequence. The difficulty lies in showing that the expansion (\ref{riccatiregexp}) carries over to our setting from the tabular MDP setting. We now turn to establishing that (\ref{riccatiregexp}) indeed is valid.

\paragraph{Riccati Equation.}
Let
\begin{align*}
x_{t+1} &= A x_t + Bu_t + w_t& x_0&=x
\end{align*}
with $x_t,w_t \in \mathbb{R}^{d_x},u_t \in \mathbb{R}^{d_u}$ and $A\in \mathbb{R}^{d_x\times d_x}, B\in \mathbb{R}^{d_x\times d_u}$ with $\E w_t w_t^\top =\Sigma_{w}\succeq 0$ and cost
\begin{align*}
V_T^\pi(x) = \E^\pi  \sum_{t=0}^{T-1}x_t^\top Q x_t + u_t^\top R u_t
\end{align*}
where $Q \succ 0, R \succ 0$. In this case, the optimal policy $\pi^*$ can be described by the following set of equations
\begin{align*}
u_t &= K_{T-t} x_t\\
K_t =&- (B^\top P_{t} B+R)^{-1} (B^\top P_{t} A)\\
P_{t+1} &=Q +A^\top P_t A - A^\top P_t B (B^\top P_t B+R)^{-1}B^\top P_t A  & P_0&=0.
\end{align*}
Above, $P_t$, by change of variables, satisfies the forward recursion (filter form). What is more, the associated optimal cost is
\begin{align*}
V_T^*(x) = x^\top P_T x \sum_{t=0}^{T-1} \tr \Sigma_{w} P_t.
\end{align*}

\paragraph{Convergence to the Stationary Cost.}
Define now
\begin{align*}
D_T(x) &= x^\top P x + \sum_{t=0}^{T-1} \E w_t P w_t = x^\top P x+ Tg & g&= \E w_t P w_t.
\end{align*}
where $P$ is stationary limit of $P_t$ and thus satisfies the algebraic Riccati equation (ARE).
\begin{align}
\label{DARE}
P &=Q +A^\top K A - A^\top PB (B^\top P B+R)^{-1}B^\top P A.
\end{align}
Direct comparison shows that
\begin{align*}
D_T - V_T^* = \sum_{t=0}^{T-1} \E w_{t} (P-P_{t+1}(T)) w_t= \sum_{t=0}^{T-1} \tr \Big[(P-P_{t+1}) \Sigma_w\Big]
\end{align*}
It follows that
\begin{align}
\label{asymptoticequivalence}
V_T = D_T +O(1)
\end{align}
where the term $O(1)$ holds uniformly in the model parameters, in a small neighborhood of $(A,B)$ preserving stability by virtue of the uniform convergence (\ref{riccaticonvergence}). Indeed, this holds uniformly on each ``strong stability'' region $\{\theta' \in \mathbb{R}^{d_\theta} | \rho(A-B(\theta')K(\theta))\leq 1 -\zeta\}$, $\zeta \in (0,1)$.

\paragraph{Rate of Convergence of the Riccati Equation.}
We state he following result concerning the rate of convergence of the Riccati recursion to its stationary limit. The following statement is adapted  from (\cite{anderson2012optimal}, p81).
\begin{theorem}
Assume that $(A,B)$ is stabilizable and that $Q\succ 0, R\succ 0$. Then
\begin{align}
\label{riccaticonvergence}
\|P_t - P\|_\infty = C [\rho(A-B(\theta)K(\theta))]^t
\end{align}
where $C>0$ can be chosen uniformly in $(A,B)$ provided the $\rho(A-BK(A,B)) <1-\zeta$ for some $\zeta \in (0,1)$.
\end{theorem}

\section{Fisher Information and Linear Dynamics}
\label{chainfisherapp}

We will need the chain rule.

\paragraph{Chain Rule for Fisher Information.}

Consider Fisher information as defined in (\ref{fisherdef}) for a bivariate density $p_\theta(x,y)$. Define the conditional Fisher information as
\begin{align*}
\I_{p(x|y)}(\theta) =\int \int \nabla_\theta \{  \log p_\theta(x|y)\} \left[\nabla_\theta\{  \log p_\theta(x|y)\} \right]^\top p_\theta(x|y)dx p_\theta(y)  dy
\end{align*}
 Then
\begin{multline}
\label{fichain}
\I_{p(x,y)}(\theta)= \int \int \nabla_\theta \log p_\theta(x,y)\left[\nabla_\theta \log p_\theta(x,y)\right]^\top p_\theta(x,y) dxdy \\
 = \int \int \nabla_\theta \log \{ p_\theta(x|y) p_\theta(y)\} \left[\nabla_\theta \log \{p_\theta(x|y) p_\theta(y)\} \right]^\top p_\theta(x|y) p_\theta(y)  dxdy\\
 =\int \int \nabla_\theta \{  \log p_\theta(x|y)+\log p_\theta(y)\} \left[\nabla_\theta \{  \log p_\theta(x|y)+\log p_\theta(y)\} \right]^\top p_\theta(x|y) p_\theta(y)  dxdy\\
 =\int \int \nabla_\theta \{  \log p_\theta(x|y)\} \left[\nabla_\theta\{  \log p_\theta(x|y)\} \right]^\top p_\theta(x|y)dx p_\theta(y)  dy\\
 +\int \int \nabla_\theta \{  \log p_\theta(y)\} \left[\nabla_\theta  \{  \log p_\theta(y)\} \right]^\top   p_\theta(x|y)   p_\theta(y)dydx\\
 +\textnormal{mixed terms linear in the expectation of score functions}\\
 =\I_{p(x|y)}(\theta) + \I_{p(y)}(\theta)
\end{multline}
assuming that the scores of $p(x|y)$ and $p(y)$ have mean zero.

\paragraph{Fisher Information and Linear Dynamics.} The chain rule presented above is useful for computing the Fisher information in the model (\ref{themodel}) of the observations $(x^{T+1}, u^T)$. We may iterate (\ref{fichain}) to find
\begin{multline}
\label{lqrfisherchain}
\I^T(\theta) = \I_{p(x^{T+1},u^T)}(\theta) = \sum_{t=0}^{T}  \E \I_t(\theta) \textnormal{ where}\\
\I_t(\theta) = \int \nabla_\theta \log p(y-A(\theta)x_t-B(\theta)u_t)\left[\nabla_\theta \log p(y-A(\theta)x_t-B(\theta)u_t)\right]^\top\\
\times p(y-A(\theta)x_t-B(\theta)u_t) dy
\end{multline}
with $p$ the common density of all noise terms $w_t$ and $y$ is a dummy variable (integration over $x_{t+1}$). The second term $I_t(\theta)$ in (\ref{lqrfisherchain}) may be reduced to
\begin{multline}
 \int \nabla_\theta \log p(y-A(\theta)x_t-B(\theta)u_t)\left[\nabla_\theta \log p(y-A(\theta)x_t-B(\theta)u_t)\right]^\top p(y-A(\theta)x_t-B(\theta)u_t) dy\\
 = \int \nabla_\theta[A(\theta)x_t+B(\theta)u_t] \nabla_w \log p(w)\left[\nabla_\theta[A(\theta)x_t+B(\theta)u_t] \nabla_w \log p(w)\right]^\top p(w) dw\\
 =\nabla_\theta[A(\theta)x_t+B(\theta)u_t]  \int \nabla_w \log p(w)\left[ \nabla_w \log p(w)\right]^\top p(w) dw \left[\nabla_\theta[A(\theta)x_t+B(\theta)u_t]\right]^\top\\
 =\nabla_\theta[A(\theta)x_t+B(\theta)u_t]  \J(p) \left[\nabla_\theta[A(\theta)x_t+B(\theta)u_t]\right]^\top
\end{multline}
where the operator $\J(\cdot)$ is defined in (\ref{fisherlocdef}). In sum, we have the following key identity
\begin{align*}
 \I^n(\theta) =\E \sum_{t=0}^n  [\dop_\theta[A(\theta)x_t+B(\theta)u_t]] ^\top \J(p) \dop_\theta[A(\theta)x_t+B(\theta)u_t]
\end{align*}
Since we assumed at $A$ does not depend on $\theta$, and $\dop \VEC B(\theta) = L$, this further reduces to
\begin{multline*}
 \I^T(\theta) =\E \sum_{t=0}^T  [ (u_t^\top \otimes I_{d_x} )\dop_\theta \VEC B(\theta)] ^\top \J(p)  (u_t^\top \otimes I_{d_x} )\dop_\theta \VEC B(\theta)\\
 =\E \sum_{t=0}^T  L^\top (u_t \otimes I_{d_x} ) \J(p)  (u_t^\top \otimes I_{d_x} ) L\\
=\E \sum_{t=0}^T  L^\top (u_tu_t^\top  \otimes \J(p)  ) L.
\end{multline*}

\section{Van Trees' Inequality}

\paragraph{Gramian Cauchy-Schwarz.}
Take two random vectors $v_1, v_2 \in \mathbb{R}^n$ and suppose that $ 0\prec \E v_2v_2^\top \prec \infty$. Observe that 
\begin{align}
\label{posdefoprod}
0\preceq
\E \begin{bmatrix}
v_1 \\ v_2
\end{bmatrix}
\begin{bmatrix}
v_1^\top & v_2^\top
\end{bmatrix}
=
\begin{bmatrix}
\E v_1v_1^\top & \E v_1 v_2^\top \\
\E v_2v_1^\top &  \E v_1 v_1^\top
\end{bmatrix}.
\end{align}
Since $\E v_2 v_2^\top \succ 0$, (\ref{posdefoprod}) implies, by the Schur complent necessary condition, that $\E v_1 v_1^\top -\E v_1v_2^\top (\E v_2v_2^\top)^{-1} \E v_2 v_1^\top \succeq 0$. Equivalently,
\begin{align}
\label{gcs}
\E v_1 v_1^\top \succeq\E v_1v_2^\top (\E v_2v_2^\top)^{-1} \E v_2 v_1^\top 
\end{align}
which is known as the Gramian Cauchy-Schwarz inequality (for scalars variabes $v_1, v_2$ it reduces to Cauchy-Schwarz). The following inequality is adapted from \cite{bobrovsky1987some}. As it is central to our main result, we provide the full details below.

\paragraph{Van Trees' Inequality.}
Let $\theta \in \mathbb{R}^p$ a random parameter with prior density $\lambda$ and let $X\in \mathbb{R}^n$ be a random vector with density $p(x|\theta)$. Define
\begin{align*}
\I(\theta) &=\int  \left( \frac{\nabla_\theta p(x|\theta)}{p(x|\theta)}\right)\left( \frac{\nabla_\theta p(x|\theta)}{p(x|\theta)}\right)^\top p(x|\theta) dx, \textnormal{ and} \\
\J(\lambda)&= \int \left( \frac{\nabla_\theta \lambda(\theta)}{\lambda(\theta)}\right) \left( \frac{\nabla_\theta \lambda(\theta)}{\lambda(\theta)}\right)^\top  \lambda(\theta)  d\theta
\end{align*}
We will now prove that
\begin{align}
\label{VTineq}
\E \left[ (\hat \psi (X)-\psi(\theta))(\hat \psi (X)-\psi(\theta))^\top\right] \succeq  \E \nabla_\theta \psi(\theta) \left[ \E \I(\theta)+\J(\lambda) \right]^{-1} \E [\nabla_\theta \psi(\theta)]^\top
\end{align}
under the regularity conditions:
\begin{enumerate}
\item[R1.] $\psi$ is differentiable.
\item[R2.] $\lambda \in C^\infty_c(\mathbb{R}^p)$; the prior is smooth with compact support.
\item[R3.] The density $p(x|\theta)$ of $X$ is continuously differentiable on the domain of $\lambda$.
\item[R4.] The score has mean zero; $ \int \left( \frac{\nabla_\theta p(x|\theta)}{p(x|\theta)}\right) p(x|\theta) dx =0$.
\item[R5.] $\J(\lambda)$ is finite and $\I(\theta)$ is a continuous function of $\theta$ on the domain of $\lambda$.
\end{enumerate}

\begin{proof}
We use (\ref{gcs}) by letting $v_1 = \hat \psi (X)-\psi(\theta)$ and $v_2 =\nabla_\theta \log [p(x|\theta) \lambda(\theta)] $. We first compute
\begin{multline*}
\E v_2 v_2^\top = \E\left[ \nabla_\theta \log [p(x|\theta) \lambda(\theta)] \left(\nabla_\theta \log [p(x|\theta) \lambda(\theta)]\right)^\top\right]\\
= \E \left(\frac{\lambda(\theta)\nabla_\theta(p(x|\theta) + p(x|\theta) \nabla_\theta \lambda(\theta)}{p(x|\theta)\lambda(\theta)} \right)\left(\frac{\lambda(\theta)\nabla_\theta(p(x|\theta) + p(x|\theta) \nabla_\theta \lambda(\theta)}{p(x|\theta)\lambda(\theta)} \right)^\top\\
=\E \left( \frac{\nabla_\theta p(x|\theta)}{p(x|\theta)}\right)\left( \frac{\nabla_\theta p(x|\theta)}{p(x|\theta)}\right)^\top +\E \left( \frac{\nabla_\theta \lambda(\theta)}{\lambda(\theta)}\right) \left( \frac{\nabla_\theta \lambda(\theta)}{\lambda(\theta)}\right)^\top+C_{p\lambda}
\end{multline*}
where
\begin{align*}
C_{p\lambda} = \E \left( \frac{\nabla_\theta p(x|\theta)}{p(x|\theta)}\right)\left( \frac{\nabla_\theta \lambda(\theta)}{\lambda(\theta)}\right)^\top +\E \left( \frac{\nabla_\theta \lambda(\theta)}{\lambda(\theta)}\right)\left( \frac{\nabla_\theta p(x|\theta)}{p(x|\theta)}\right)^\top=0
\end{align*}
by the mean zero property of the score. Hence
\begin{align}
\label{v2v2eqn}
\E v_2 v_2^\top = \E \I(\theta)+\J(\lambda).
\end{align}
We still need to establish that the matrix in (\ref{v2v2eqn}) has full rank. This is actually a byproduct of the computation of $\E v_1 v_2^\top$. Now since $\lambda\in C^\infty_c(\mathbb{R}^p)$ is a test function, we may integrate by parts to find
\begin{multline}
\label{ibp}
\E v_1 v_2^\top = \E \left[ (\hat \psi (X)-\psi(\theta))\left(\nabla_\theta \log [p(x|\theta) \lambda(\theta)] \right)^\top \right]\\
=\E \left[ (\hat \psi (X)-\psi(\theta))\left(\nabla_\theta \log [p(x,\theta) ] \right)^\top \right]=\E \left[ (\hat \psi (X)-\psi(\theta))\left(\frac{\nabla_\theta p(x,\theta)}{p(x,\theta)}  \right)^\top \right]\\
=\int\int (\hat \psi (X)-\psi(\theta))\left(\frac{\nabla_\theta p(x,\theta)}{p(x,\theta)}  \right)^\top p(x,\theta) dx d\theta\\
=\int\int \nabla_\theta \psi(\theta) p(x,\theta) dx d\theta = \E \nabla_\theta \psi(\theta).
\end{multline}
In particular, using $\psi(\theta) = \theta$ yields $\E v_1v_2 =I$, which is sufficient to conclude that $\E v_2 v_2$ has full rank. The result follows by (\ref{gcs}) combined with (\ref{v2v2eqn}) and (\ref{ibp}).
\end{proof}

\section{Matrix Calculus and Spectral Pertubation}
\label{appmatcal}
\subsection{Matrix Calculus}
It is frequently required of us to compute derivatives of the form
\begin{align*}
\dop_{\theta} [M(\theta) v]
\end{align*}
where $\theta \in \mathbb{R}^{d_\theta}, M \in \mathbb{R}^{m\times n}$ and $v\in \mathbb{R}^{n}$. Observe that we are not actually taking the derivative with respect to the matrix $M(\theta) \in \mathbb{R}^{m\times n}$ but rather the vector $[M(\theta) v] \in \mathbb{R}^m$. However, since $v$ does not depend on $\theta$, one expects to be able to separate the variables $M$ and $v$.  Some facts about this operation are collected here.

\paragraph{Useful Identities.}

Let $M, N , P $ be three matrices such that product MNP exists. A useful formula which we shall make use of below is
\begin{align}
\label{vecid}
\VEC MNP = (P^\top \otimes M)\VEC N
\end{align}

\paragraph{Derivative of Matrix-Vector Product.} Observe that by the identity (\ref{vecid}) just stated that
\begin{align}
\VEC M(\theta) v = (v^\top \otimes I_m) \VEC M(\theta).
\end{align}
Hence 
\begin{align}
\label{dopMvformula}
\dop_{\theta} [M(\theta) v] =   (v^\top \otimes I_m)[\dop_{\theta} \VEC M(\theta)] .
\end{align}
We also require the following result.

\begin{lemma}
\label{kronprodlemma}
Let $M(\theta)$ be a smooth matrix function and $N,P \succ 0$ be matrices of appropriate size. Then for any vector $v$
\begin{align*}
\tr \left( N [\dop_\theta M(\theta) v] P  [\dop_\theta M(\theta) v]^\top \right)=\tr \left([( vv^\top \otimes N)   [\dop_\theta \VEC M(\theta) ]P   \dop_\theta \VEC  M(\theta) ]^\top \right).
\end{align*}
\end{lemma}
\begin{proof}
Expand $0 \prec N = \sqrt{N} \sqrt{N}$ and use identity (\ref{vecid}) in combination with (repeated use of) the trace cyclic property to find
\begin{multline*}
\begin{aligned}
&\tr \left( N [\dop_\theta M(\theta) v] P  [\dop_\theta M(\theta) v]^\top \right)\\&=\tr \left( \sqrt{N} \sqrt{N} [\dop_\theta M(\theta) v]P  [\dop_\theta M(\theta) v]^\top \right)\\
&=\tr \left(   [\dop_\theta\sqrt{N} M(\theta) v]P  [\dop_\theta\sqrt{N}M(\theta) v]^\top  \right)\\
&=\tr \left(   [( v^\top \otimes \sqrt{N})\dop_\theta \VEC M(\theta) ]P   [( v^\top \otimes \sqrt{N})\dop_\theta \VEC  M(\theta) ]^\top  \right)\\
&=\tr \left([\dop_\theta \VEC  M(\theta) ]^\top( v \otimes \sqrt{N})   [( v^\top \otimes \sqrt{N})\dop_\theta \VEC M(\theta) ]P    \right)\\
&=\tr \left([\dop_\theta \VEC  M(\theta) ]^\top( vv^\top \otimes N)   [\dop_\theta \VEC M(\theta) ]P    \right)\\
&=\tr \left([( vv^\top \otimes N)   [\dop_\theta \VEC M(\theta) ]P   \dop_\theta \VEC  M(\theta) ]^\top \right)
\end{aligned}
\end{multline*}
as desired.
\end{proof}

\paragraph{Linear Parametrizations.} Assume that the map $\theta \mapsto M(\theta)$ is affine. In this case the derivative of $M(\theta)$ is entirely characterized by a linear map $L\in \mathbb{R}^{mn \times d_\theta}$ such that
\begin{align*}
\VEC M(\theta) = L\theta.
\end{align*}
In this case the formula (\ref{dopMvformula}) becomes particularly simple, and we see that it reduces to
\begin{align*}
\dop_{\theta} [M(\theta) v] =   (v^\top \otimes I_m)L \in \mathbb{R}^{n\times d_\theta} .
\end{align*}
A useful reference for these results is \cite{magnus2019matrix}.

\subsection{Proof of Proposition~\ref{charprop}}
\label{charpropproof}
In particular, the above tools readily yield the following proposition.
\begin{proof}
Write, for each $t$, 
\begin{multline}
[\dop_\theta[B(\theta)u_t]]^\top  \J(p) \dop_\theta[B(\theta)u_t] = [ (u_t^\top \otimes I_{d_x})\dop_{\theta} \VEC B(\theta)]^\top  (u_t^\top \otimes \J(p))\dop_{\theta} \VEC B(\theta)\\
 =[ \dop_{\theta} \VEC B(\theta) ]^\top (u_t \otimes I_{d_x})  (u_t^\top \otimes \J(p))\dop_{\theta} \VEC B(\theta)\\
 =[ \dop_{\theta} \VEC B(\theta) ]^\top [u_tu_t^\top \otimes  \J(p) ]\dop_{\theta} \VEC B(\theta).
\end{multline}
Since $u_t = Kx_t$ and $\ker \E Kx_tx_t^\top K^\top = \ker KK^\top$ using $\E x_t x_t^\top \succeq \Sigma_w\succ 0$ the result is established.
\end{proof}

\subsection{Proof of the claim in Example~\ref{memorylessexample}}
We first establish the following auxilliary result, which studies the infinitesimal perturbations yielding the same optimal policy in terms of an eigenvalue equation.

\begin{lemma}
\label{ch4:kerklemma}
Under the hypotheses of Example~\ref{memorylessexample}, the nullspace $\ker  \dop_\theta \VEC K(\theta)$ is given by the eigenspace associated to the eigenvalue $1$ of
\begin{align}
\label{ch4:goodsamaritan}
\left[ ( K^\top F^\top \otimes I_{d_u})\Pi_{d_u d_y}   +(K^\top \otimes F^\top) \right].
\end{align}
where the permutation matrix $\Pi_{d_u d_y}$ maps  $ \VEC  M$ to $\VEC M^\top$ for $M \in \mathbb{R}^{d_u\times d_y}$.
\end{lemma}

\begin{proof}
We compute the differential of $K$ with respect to $F$ as
\begin{multline*}
\diff (F^\top F +\lambda I_{du})^{-1} F^\top\\
= (F^\top F+ \lambda I_{d_u})^{-1}(\diff F)^\top - (F^\top F+\lambda I_{d_u})^{-1}[(\diff F)^\top F + F^\top \diff F] (F^\top F+\lambda I_{d_u})^{-1}F^\top.
\end{multline*}
Multiplying by $(F^\top F+\lambda I_{d_u})$, we see that this can be set to zero if and only if there is a root $dF$ to the equation
\begin{align*}
(\diff F)^\top = [(\diff F)^\top F + F^\top \diff F] (F^\top F+\lambda I_{d_u})^{-1}F^\top.
\end{align*}
Vectorizing using identity (\ref{vecid}), we find the equivalent equation
\begin{multline*}
\VEC (\diff F)^\top = \VEC \left( [(\diff F)^\top F + F^\top \diff F] (F^\top F+\lambda I_{d_u})^{-1}F^\top\right)\\
=( F  (F^\top F+ \lambda I_{d_u})^{-1} F^\top \otimes I_{d_u}) \VEC \diff F^\top +(F(F^\top F+\lambda I_{d_u})^{-1} \otimes F^\top) \VEC \diff F.
\end{multline*}
Introducing $\Pi$ such that $\VEC (\diff F)^\top = \Pi_{d_u d_y}  \VEC \diff F$, this is in turn equivalent to the eigenvalue equation
\begin{multline*}
\begin{aligned}
&\VEC \diff F \\
=& \Pi_{d_u d_y} ^{-1} \left[ ( F  (F^\top F+ \lambda I_{d_u})^{-1} F^\top \otimes I_{d_u})\Pi_{d_u d_y}  \VEC \diff F +(F(F^\top F+\lambda I_{d_u})^{-1} \otimes F^\top) \VEC \diff F\right]\\
=&  \Pi_{d_u d_y} ^{-1} \left[ ( K^\top F^\top \otimes I_{d_u})\Pi_{d_u d_y}  \VEC \diff F +(K^\top \otimes F^\top) \VEC \diff F\right].
\end{aligned}
\end{multline*}
This completes the proof.
\end{proof}

\begin{proposition}
\label{ch4:propfromex}
The problem given in Example~\ref{memorylessexample} is uninformative if and only if $KK^\top $ is singular. Moreover the information singular subspace $\mathtt{U}$ is unique and has dimension
\begin{align*}
\dim \mathtt{U} =  d_y \times \dim  \ker KK^\top .
\end{align*}
\end{proposition}

\begin{proof}
Suppose first that the instance is uninformative. Now, since $\J(p)$ has full rank, $KK^\top$ must necessarily be singular for $KK^\top \otimes \J(p)$ to be singular.

For the other direction, we will prove that 
\begin{align}
\label{ch4:ucand}
\mathtt{U}_{cand} = \{ \tilde u \otimes \tilde w \in \mathbb{R}^{d_\theta} : \tilde u \in \ker KK^\top, \tilde w \in \mathbb{R}^{d_y} \} 
\end{align}
is an information singular subspace. It is clear that $\mathtt{U}_{cand} = \ker KK^\top \otimes \J(p)$, so every other candidate subspace would be a proper subspace of $\mathtt{U}_{cand}$, rendering it unique. To finish the proof, it thus suffices to prove that no nonzero vector of $\mathtt{U}_{cand}$ is in $\ker  \dop_\theta \VEC K(\theta)$. The idea is to use Lemma~\ref{ch4:kerklemma}, for which we must establish that no nonzero vector $\tilde v \in \mathtt{U}_{cand}$  belongs to the eigenspace of the eigenvalue $1$ of (\ref{ch4:goodsamaritan}).

To that end,  consider any vector of the form $\tilde v = \tilde u \otimes \tilde w $ where $\tilde u \in \ker KK^\top = \ker K^\top$ and $\tilde w \in \mathbb{R}^{d_y}$. Obviously $\tilde v \in \ker  K K^\top  \otimes  \J(p)$. Moreover, the action of $\Pi_{d_u d_y}$ on vectors $\tilde u \otimes \tilde w$ is described by $\Pi_{d_u d_y}\tilde u \otimes \tilde w = \tilde w\otimes \tilde u$. We thus compute
\begin{equation}
\begin{aligned}
\label{ch4:anothersamaritan}
&\left[ (K^\top F^\top \otimes I_{d_u})\Pi_{d_u d_y}   +(K^\top \otimes F^\top) \right] (\tilde u\otimes \tilde w)\\
 &= ( K^\top F^\top \otimes I_{d_u})\Pi_{d_u d_y} (\tilde u\otimes \tilde w)  + (K^\top \otimes F^\top) (\tilde u\otimes \tilde w) \\
&=( K^\top F^\top \otimes I_{d_u})(\tilde w\otimes \tilde u)  + (K^\top \otimes F^\top) (\tilde u\otimes \tilde w)\\
&=( K^\top F^\top \tilde w\otimes I_{d_u} \tilde u) + (\underbrace{K^\top\tilde u}_{=0} \otimes F^\top\tilde w)\\
&=( K^\top F^\top \otimes I_{d_u}) (\tilde w \otimes  \tilde u).
\end{aligned}
\end{equation}
using $\ker KK^\top = \ker K^\top$.

We shall now prove that the spectrum of $( K^\top F^\top \otimes I_{d_u}) $ is contained in $(-\infty, 1)$. Since this does not include $1$, it then follows by Lemma~\ref{ch4:kerklemma} that all nonzero elements of the subspace (\ref{ch4:ucand}) are not in $\ker  \dop_\theta \VEC K(\theta)$. To show this, it suffices to prove that  the spectrum of $( K^\top F^\top) $ is contained in $(-\infty, 1)$ due to the structure of the eigenvalues for Kroenecker products. For the next step, notice that since any $\tilde w \in \ker K^\top F^\top$ by (\ref{ch4:anothersamaritan}) belongs to the nullspace of (\ref{ch4:goodsamaritan}), we may restrict attention to $\tilde w \in (\ker K^\top F^\top)^\perp$. Notice now that for all such $\tilde w$ and for $\lambda' \in (0,\lambda)$, we have that
\begin{align}
\label{ch4:monorelation}
\tilde w^\top K^\top F^\top\tilde w =\tilde w^\top  F  (F^\top F +\lambda I_{d_u})^{-1} F^\top \tilde w < \tilde w^\top  F  (F^\top F +\lambda' I_{d_u})^{-1} F^\top\tilde w.
\end{align}
However, 
\begin{align*}
\lim_{\lambda'\to 0}  F  (F^\top F +\lambda' I_{d_u})^{-1} F^\top \tilde w =F^\dagger F \tilde w = \tilde w
\end{align*}
since $F^\dagger F$ is the orthogonal projector on the range of $F$ and $\tilde w  \in (\ker K^\top F^\top)^\perp =\im  FK \subset \im F$, by range-nullspace. Combining this observation with (\ref{ch4:monorelation}) we obtain the inequality
\begin{align*}
\tilde w^\top K^\top F^\top\tilde w < \tilde w^\top \tilde w,
\end{align*}
valid for all $\tilde w \in (\ker K^\top F^\top)^\perp$, which can only be true if the largest eigenvalue of $K^\top F^\top$ is smaller than $1$ (note that the eigenvalues are nonnegative real since $ K^\top F^\top \succeq 0$). 

It follows by Lemma~\ref{ch4:kerklemma} that every nonzero $\tilde w \in \mathtt{U}_{cand}$ satisfies (\ref{uninformativechar}), and it is easly seen that $\dim \mathtt{U}_{cand} =  d_y \times \dim  \ker KK^\top$, which is the dimension of the space $\ker  [K(\theta) K(\theta)^\top  \otimes  \J(p) ]$. Comparing with  (\ref{uninformativechar}), $\dim \mathtt{U}_{cand}$ has maximal dimension and so must be unique (all other candidates would be proper subspaces of  $\ker KK^\top \otimes \J(p)$).
\end{proof}

\subsection{Davis-Kahan $\sin \theta$ Theorem}
Consider any symmetric positive semidefinite matrix $N \in \mathbb{R}^{n\times n}$. Its spectral decomposition may be written as
\begin{align}
\label{dkdec1}
N = V_1 \Lambda_1 V_1^\top + V_0 \Lambda_0 V_0^\top
\end{align}
where $\Lambda_1= \diag(\lambda_1,\dots, \lambda_{I})$, $\Lambda_0=\diag(\lambda_{I+1},\dots, \lambda_n)$ and $V_1, V_0$ are partial isometries with columns spanning the corresponding eigenspaces (the obvious choice here is $V_i=O_i$, with $O_i$ having orthonormal columns). Consider now the matrix $\tilde N = N+T$, where $T$ is a ``small'' symmetric perturbation. Clearly, we may write
\begin{align}
\label{dkdec2}
 \tilde N =  \tilde V_1 \tilde \Lambda_1\tilde V_1^\top + \tilde V_0 \tilde \Lambda_0  \tilde V_0^\top
\end{align}
with the variables $\tilde V_i,\tilde \Lambda_i$ defined analogously.

The Davis-Kahan $\sin\theta$-Theorem concerns itself with controlling the deviations between $V_i$ and $\tilde V_i$ in terms of the magnitude of the pertubation $T$. Define for any unitarily invariant norm $\|\cdot\|$ and any two subspaces, $S,\tilde S$
\begin{align*}
\| \sin \theta (S,\tilde S)\|=\|(I-\pi_S)\pi_{\tilde S}\|
\end{align*}
where $\pi_S, \pi_{\tilde S}$ are the orthogonal projections onto the subspaces $S,\tilde S$. Note that this definition is symmetric in $S,\tilde S$ since orthogonal projections commute. We extend this definition to any two matrices $M,\tilde M$ by 
\begin{align*}
\| \sin \theta (M,\tilde M)\|=\|(I-\pi_{\spn M})\pi_{\spn\tilde M}\|
\end{align*}
We now state without proof a version of the Davis-Kahan $\sin\theta$ theorem most amenable to our needs \citep{wedin1972perturbation}.

\begin{theorem}
Let $N$ and $\tilde N$ by symmetric positive semidefinite matrices with spectral decompositions (\ref{dkdec1}) and (\ref{dkdec2}) respectively. Assume there exist $\alpha\geq 0$ and $\delta >0$ such that $\sigma_{\min}(\tilde \Lambda_1)\geq \alpha+\delta$ and $\sigma_{\max}(\Lambda_0)\leq \alpha$. Then for every unitarily invariant norm
\begin{align*}
\|\sin\theta(V_0,\tilde V_0)\| \leq \frac{2\|T\tilde V_0  \|}{\delta}.
\end{align*}
\end{theorem}
more statements are possible, but this version is sufficient for our purposes. An up-to-constants-equivalent formulation of the $\sin\theta$ distance for the spectral norm $\|\cdot\|_\infty$ is
\begin{align}
\label{ddef}
d_{\infty}(V,\tilde V) =\inf_{O\in\mathbb{O}_{n} } \|V -\tilde V O\|_\infty 
\end{align}
for any $V,\tilde V \in \mathbb{O}_{m,n}$ (real matrices with orthonormal columns). We have the following equivalence statement \citep{cai2018rate}.
\begin{proposition}
For any two $V,\tilde V \in\mathbb{O}_{m,n}$, it holds that
\begin{align*}
\|\sin\theta(V,\tilde V)\|_\infty \leq  d_{\infty}(V,\tilde V) \leq \sqrt{2} \|\sin\theta(V,\tilde V)\|_\infty.
\end{align*}
\end{proposition}

%\vskip 0.2in

%\bibliography{sample}

\end{document}